\documentclass[a4paper,10pt]{article}
\usepackage{romannum}
\usepackage{subcaption}
\usepackage{epstopdf}
\usepackage{graphics,graphicx}
\usepackage{amsmath,amsthm,amssymb}

\title{An oversampled collocation approach of the Wave Based Method for Helmholtz problems}

\author{Daan Huybrechs and Anda-Elena Olteanu\footnote{KU Leuven, Department of Computer Science, Celestijnenlaan 200A, BE-3001 Leuven, Belgium. Corresponding author: daan.huybrechs@cs.kuleuven.be}}

\begin{document}

\pagenumbering{arabic}

\maketitle

\begin{abstract}
The Wave Based Method (WBM) is a Trefftz method for the simulation of wave problems in vibroacoustics. Like other Trefftz methods, it employs a non-standard discretisation basis consisting of solutions of the partial differential equation (PDE) at hand. We analyse the convergence and numerical stability of the Wave Based Method for Helmholtz problems using tools from approximation theory. We show that the set of discretisation functions more closely resembles a frame, a redundant set of functions, than a basis. The redundancy of a frame typically leads to ill-conditioning, which indeed is common in Trefftz methods. Recent theoretical results on frames for function approximation suggest that the associated ill-conditioned system matrix can be successfully regularised, with error bounds available, when using a discrete least squares approach. While the original Wave Based Method is based on a weighted residual formulation,  in this paper we pursue an oversampled collocation approach instead. We show that, for smooth scattering obstacles in two dimensions, the results closely follow the theory of frames. We identify cases where the method achieves very high accuracy whilst providing a solution with small norm coefficients, in spite of ill-conditioning. Moreover, the accurate results are reliably maintained even in parameter regimes associated with extremely high ill-conditioning.
\end{abstract}

\section{Introduction}
\label{sect:introduction}

Numerical simulation methods for wave scattering and propagation problems lead to a wealth of mathematical and computational challenges. One of the primary concerns is the number of degrees of freedom that are required in order to represent a wave field, in particular for problems involving high or moderately high frequencies. The class of Trefftz methods aims to reduce this number, by discretising the governing partial differential equation (PDE) using solutions of that same PDE. We refer the reader to \cite{hiptmair2016trefftz} and references therein for a review of Trefftz methods, such as the Method of Fundamental Solutions \cite{fairweather2003mfs}, the Wave Based Method \cite{desmet98phd,desmet2014wbm} and the Ultra-Weak Variational Formulation scheme \cite{cessenat1998uwvf}. Since the ansatz of the solution satisfies the PDE by construction, it remains to enforce boundary and interface conditions. This leads to a linear system of equations that is typically dense, albeit with much smaller dimensions than a corresponding discretisation using finite elements. Unfortunately, in quite a few cases, the system exhibits a large or very large condition number~\cite{hiptmair2016trefftz}. This seems a cause of concern, and perhaps even prevents wider adoption of these methods. Yet, extensive experiments with and literature on Trefftz methods seem to indicate that, in spite of ill-conditioning of the linear system, high accuracy solutions are often found and high-order convergence is observed.

The analysis of Trefftz methods has seen great progress in recent years. A major focus of the mathematical literature on this topic lies with the best approximation to the solution in the approximation space, the space that is spanned by the chosen basis functions \cite{hiptmair2014harmonic,hiptmair2016pwdg}. Given the unconventional nature of the space, this is a non-trivial problem. The computation of the approximation itself, in particular the way it is obtained from a potentially ill-conditioned system, has received much less attention. Barnett and Betcke have studied the Method of Fundamental Solutions, and proposed an implementation strategy that is numerically stable by a judicious choice of the charge points in the method \cite{barnett2008stability,barnett2010nonpolynomial}. An often quoted test is the Picard condition \cite[\S1.2.3]{hansen2005illposed}. This is a test that can provide confidence in a solution a-posteriori, but does not provide an a-priori guarantee of success.

Some recent methods in approximation theory have similar characteristics as Trefftz methods, similar in the sense that highly accurate solutions are found by solving ill-conditioned systems. These methods are based on using a so-called \emph{frame}, rather than a basis. We review frames in \S\ref{sect:frames} of this paper. Compared to a basis a frame is redundant, hence there are multiple representations of any given function in the frame. This leads to ill-conditioned matrices for the approximation problem. Yet, generic error bounds can be shown \cite{huybrechs2016tw674,adcock2018frames2}. Moreover, convergence to high accuracy, up to machine precision even, can provably be guaranteed in spite of the ill-conditioning, under certain additional conditions. The main condition is \emph{oversampling}: rather than solving square systems, accuracy and robustness is significantly improved by considering solutions in a least squares sense. Best results are achieved in particular with a discrete least squares approximation, based on a larger number $M$ of function evaluations than degrees of freedom $N$ (i.e., $M > N$, leading to a rectangular system). This approach is shown to be successful precisely when the approximation space is obtained from a frame. Hence the importance of the concept of frames in this context: if the redundant set is not a frame, functions can be found for which the approximation scheme may not converge using finite precision computations. 

The analogue of a discrete least-squares approximation for Trefftz methods is an oversampled collocation approach. The goal of this paper is to describe such a collocation approach for the Wave Based Method, for Helmholtz problems in two dimensions. We fully address the problem of ill-conditioning of the linear systems and, through a sequence of experiments, illustrate how it relates to the error bounds from the literature on frames. The experiments are inspired by the prior analysis of the Method of Fundamental Solutions in \cite{barnett2008stability} and, following Barnett and Betcke, we interpret the requirement of having a convex domain in WBM in terms of external singularities of the solution. We show much improved accuracy of the oversampled collocation approach, compared to the weighted residual formulation that is commonly employed. In addition, the approach is more efficient since no calculation of integrals is required in the setup of the system matrix.

On the other hand, our results are restricted to convex and certain non-convex domains with a smooth boundary in two dimensions. Whether the beneficial properties of an oversampled collocation approach remain in the various other settings in which WBM has been applied, in particular in the presence of corner singularities, is not certain and is a topic of further research.

We also like to stress that, though we include no mathematical analysis in this paper, the set of functions used in WBM does \emph{not} satisfy the mathematical conditions of a frame. The means that, mathematically, an accurate numerical solution is not guaranteed for all possible boundary conditions. Correspondingly, we show examples where full accuracy is not obtained.

Ongoing research efforts seem to indicate that this is the case for several other Trefftz methods as well, even when it is known that the best approximation from the approximation space converges. Yet, the connection to the literature on frames does support the following observation. It is not sufficient that a good approximation to the solution exists in the approximation space. However, it \emph{is} sufficient that a good approximation exists that has expansion coefficients with (sufficiently) small discrete norm. If a good approximation with moderately small coefficients exists then, with sufficient oversampling and using a sufficiently large number of degrees of freedom, a numerical solution will be found with equal or with better accuracy. Perhaps surprisingly, this is true regardless of the ill-conditioning of the linear system. We believe this observation may underlie the success of Trefftz methods in a variety of practical applications. In this paper we support the statements above with numerical experiments and interpretations for WBM.

The structure of the paper is as follows. In \S\ref{sect:wbm} we recall a formulation of the Wave Based Method. We review some relevant concepts of the theory of frames in \S\ref{sect:frames}. We describe in particular the generic error bounds on which our interpretations are based. We formulate an oversampled collocation implementation of the Wave Based Method in \S\ref{sect:collocation}. We present a sequence of numerical experiments and their results in \S\ref{sect:results}. The results are interpreted using the above-mentioned error bounds in \S\ref{sect:discussion}. Finally, we make some concluding remarks in \S\ref{sect:conclusions}.

\section{The Wave Based Method}
\label{sect:wbm}
 
We consider the homogeneous Helmholtz equation in the interior of a smooth and bounded domain $\Omega\subset\mathbb{R}^2$,
 \begin{equation}\label{eq:Helmholtz}
 \Delta u(x,y)+ k^2u(x,y)=0, \quad (x,y) \in \Omega,
 \end{equation}
 with a given real constant wavenumber $k$, $[\hbox{m}^{-1}]$, and a corresponding boundary condition, e.g. of the commonly used Neumann type
  \begin{equation}\label{eq:boundary}
  \displaystyle\frac{\partial u}{\partial n}= w(s),\quad s\in\partial\Omega.
  \end{equation}
The meaning of the unknown function $u$ depends on the wave phenomenon modeled at hand, and typically in the acoustic field, where the WBM has been developed, it denotes the steady state acoustic pressure.

Since the WBM belongs to the family of Trefftz methods, the solution of equation \eqref{eq:Helmholtz} on the domain of interest $\Omega$ is sought as an expansion of a complete set of wave functions, which satisfy exactly the governing equation \cite{desmet98phd,desmet2014wbm}.
 
\begin{figure}[tbp]
   	\centering
   	\includegraphics[scale = 0.7]{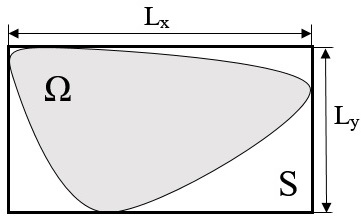}
   	\caption{Circumscribing box $S$ enclosing the domain $\Omega.$}
   	\label{fig:WBMdomains}
   \end{figure}

In the case of WBM for Helmholtz problems in two dimensions, the functions are denoted by $\{\Phi^{(\pm)}_{m},\Phi^{(\pm)}_{n}\}_{m,n=0}^{\infty}$. They are defined naturally on a rectangular domain $S\subset\mathbb{R}^2$ that circumscribes the domain of interest $\Omega$. We assume without loss of generality that the box can be chosen to be $S = [0,L_x]\times[0,L_y]$, see Fig.~\ref{fig:WBMdomains}.
Thus, the solution $u(x,y)$ is approximated by $\tilde{u}(x,y)$ given by
 \begin{equation}\label{eq:wbmSol}
 u\approx \tilde{u} =\sum\limits_{m=0}^{N_m}\left(\alpha_m^{(1)}\Phi_m^{(+)} + \alpha_m^{(2)}\Phi_m^{(-)}\right)+     \sum\limits_{n=0}^{N_n}\left(\alpha_n^{(3)}\Phi_n^{(+)} + \alpha_n^{(4)}\Phi_n^{(-)}\right),
 \end{equation}
where $N_m$ and $N_n$ represent the number of functions used in the $x$ and $y$ directions, respectively. The wave functions themselves are given explicitly as
\renewcommand{\arraystretch}{1.8}
\begin{equation}\label{eq:wavef}
\begin{array}{l l}
\Phi_m^{(\pm)} = \cos\left(k_x^{(m)}x\right)\hbox{e}^{-ik_y^{(m)}y},\\
\Phi_n^{(\pm)} = \hbox{e}^{-ik_x^{(n)}x}\cos\left(k_y^{(n)}y\right),\\
\end{array}
\end{equation}
with the wavenumber components
\begin{equation}\label{eq:wavenumb}
\begin{array}{l}
\left(k_x^{(m)},k_y^{(m)}\right)= \left(\displaystyle\frac{m\pi}{L_x},\,\,\pm\sqrt{k^2-\left(\frac{m\pi}{L_x}\right)^2}\right), \qquad m=0,\ldots,N_m,\\
\left(k_x^{(n)},k_y^{(n)}\right)= \left(\pm\sqrt{k^2-\left(\displaystyle\frac{n\pi}{L_y}\right)^2},\,\,\displaystyle\frac{n\pi}{L_y}\right), \qquad n=0,\ldots,N_n,\\
\end{array}
\end{equation}
such that $\left({k_x^{(m)}}\right)^2 + \left({k_y^{(m)}}\right)^2 = \left({k_x^{(n)}}\right)^2 + \left({k_y^{(n)}}\right)^2 = k^2$ holds and \eqref{eq:Helmholtz} is satisfied by the expansion functions. These are standing waves in one direction, and either propagating or evanescent waves in the other direction, depending on whether the corresponding wavenumber component is real or purely imaginary.

It should be noted that a complex value of the wavenumber components $k_y^{(m)}$ or $k_x^{(v)}$ of the wave functions $\Phi_m^{(+)}$ and $\Phi_n^{(+)}$ may lead to amplitudes significantly larger or significantly smaller than $1$. For numerical considerations, these functions require scaling \cite{desmet98phd}. Regarding truncation, 
it is a rule of thumb to choose a user defined truncation parameter $T$ such that the wavenumber components are linearly scaled with the physical wavenumber $k$ by the factor $T$, leading to
\begin{equation}\label{eq:truncation}
N_m = \left\lceil\frac{kTL_x}{\pi}\right\rceil \quad \mbox{and} \quad N_n = \left\lceil\frac{kTL_y}{\pi}\right\rceil.
\end{equation}
 
The problem then translates to finding the unknown coefficients of the expansion, enumerated as a single vector $\alpha = [\alpha_m^{(1)}\,\,\alpha_m^{(2)}$ $\,\,\alpha_n^{(3)}\,\,\alpha_n^{(4)}]$, $m = 0,\ldots,N_m,$ $n = 0,\ldots,N_n$. We enumerate the basis functions accordingly. A solution is found by imposing the boundary condition in a weighted residual formulation:
\begin{equation}\label{eq:weightRes}
\int\limits_{\partial\Omega}\Phi\left(\displaystyle\frac{\partial u}{\partial n}-w\right)ds=0.
\end{equation}
This leads to a linear system of equations $A\alpha = b$ of size $N=2(N_m+1)+ 2(N_n+1),$ where 
\begin{equation}\label{eq:integrals}
(a_{i,j})=\int\limits_{\partial\Omega}\Phi_i\displaystyle\frac{\partial\Phi_j}{\partial n}\hbox{d}s,\quad(b_{i})=\int\limits_{\partial\Omega}\Phi_iw\hbox{d}s,\quad i,\,\,j = 1,\ldots, N.
\end{equation}
The resulting matrix $A$ is symmetric, fully populated, complex valued and increasingly ill-conditioned as more wave functions are added to the approximation. Moreover, since the system is set up through evaluation of integrals, care needs to be taken to ensure sufficient accuracy in computations.
 
Regarding convergence, a sufficient condition stated in \cite{desmet98phd} for the approximation $\tilde{u}$ to converge to the exact solution $u$ as $N\rightarrow\infty$ is that the domain $\Omega$ is convex. When this is not the case, several alternative approaches are available. A standard practice is to decompose the given domain into several convex subdomains, or, for more complex geometries, one can employ a multi-level approach \cite{genechten2010multilevel}, or hybrid variants of the WBM \cite{vanHal04phd}. 

\section{Frames and numerical function approximation}
\label{sect:frames}

The accuracy of a numerical simulation scheme for PDE's is greatly influenced by the extent to which the solution can be approximated by a linear combination of the chosen basis functions. In particular, the numerical solution can not have greater accuracy than the best approximation to the solution in the approximation space, which is the span of the basis functions. Since Trefftz methods, and WBM in particular, employ non-standard basis functions (different from piecewise polynomials, say, which are well understood), we introduce a concept from approximation theory, \emph{frames}. This will allow to investigate and explain the convergence behaviour of WBM.

\subsection{Definition and main properties}
A frame $\Phi := \{ \phi_k \}_{k=1}^\infty$ is a set of functions that is complete in a function space, such as the space $L^2(\Omega)$ of square integrable functions on a domain $\Omega$. In addition, it should satisfy the so-called \emph{frame property}:
\begin{equation}\label{eq:frame}
 A \Vert f \Vert_{L^2(\Omega)}^2 \leq \sum_{k=1}^\infty \langle f, \phi_k \rangle^2 \leq B \Vert f \Vert_{L^2(\Omega)}^2, \qquad \forall f \in L^2(\Omega).
\end{equation}
This condition seems technical and not immediately intuitive: the main consequence of the frame property is that each function $f$ in the space can be written as a linear combination of the frame elements $\phi_k$, with coefficients that are bounded in norm. The latter restriction, which we return to several times in this paper, is a crucial difference between a frame and any other redundant complete set of functions. We refer to \cite{christensen2008framesandbases} for a full mathematical treatment of frames, and to \cite{huybrechs2016tw674} for a recent review of the use of frames for numerical function approximation. In what follows, we focus on the main conclusions in \cite{huybrechs2016tw674}.

\begin{figure}
\begin{center}
 \includegraphics[width=5cm]{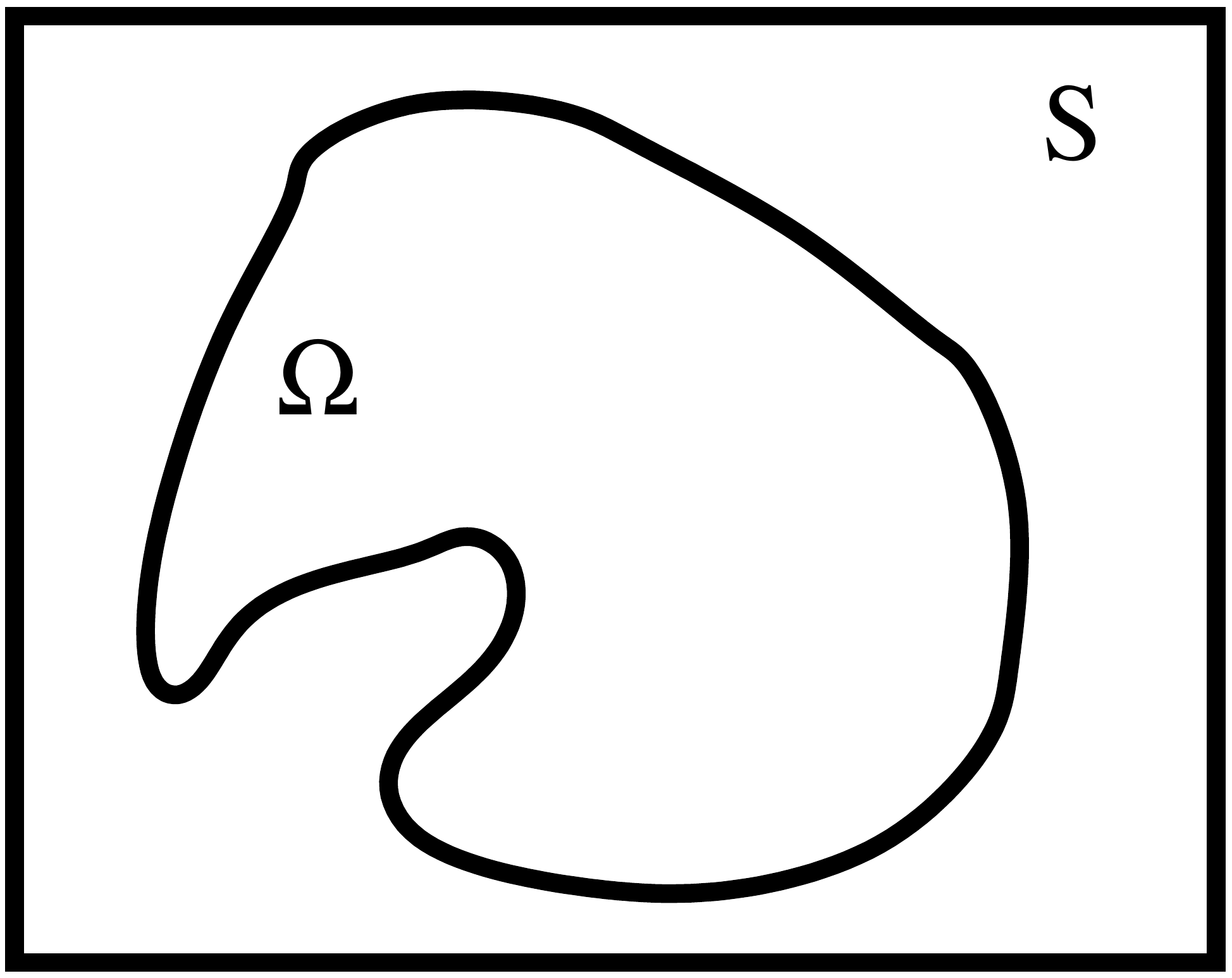}
 \caption{If $\Omega \subset S$ is a subdomain of the box $S$, then any basis for $L^2(S)$ is a frame for $L^2(\Omega)$~\cite{huybrechs2016tw674}. This is true regardless of the geometric complexity of $\Omega$. A popular choice on the box $S$ is a tensor-product Fourier series. The approximation of a given function on $\Omega$ by a Fourier series on a larger bounding box is called the \emph{Fourier extension problem}. While the construction of a basis on $\Omega$ is difficult for complex shapes, often leading to meshes and low-order piecewise polynomials, the construction of a frame is nearly trivial.}\label{fig:extensionframe}
\end{center}
\end{figure}

Of importance for the current paper is the fact that frames may be overcomplete or redundant. Indeed, condition \eqref{eq:frame} does not even prevent linear dependencies between the functions $\phi_k$. A simple frame with redundancy is obtained by adding two bases for the same space. Or, a basis can be augmented with a small number of additional functions. In both cases, a single basis is already complete and any addition to it is redundant. A third example is illustrated in Fig. \ref{fig:extensionframe}: a frame on a domain can always be obtained from a basis on a larger domain. This setting is mathematically not the same as the construction of the basis functions in WBM, but it is certainly related.

Frames are much more flexible than bases, as the example of Fig.~\ref{fig:extensionframe} shows where the domain $\Omega$ can have any shape. Unfortunately, the computation of the approximation to a function with $N$ frame elements requires solving a linear system
\[
Gc=B,
\]
where the solution vector $c$ yields the expansion coefficients for an approximation of the form $f \approx \sum_{k=1}^N c_k \phi_k(x)$. That is, unlike for an orthogonal basis, the coefficients $c_k$ are not simply given by the inner products $\langle f, \phi_k \rangle$ (even though those inner products do appear in the frame property \eqref{eq:frame}). Moreover, if the frame is redundant, the matrix $G$ is highly ill-conditioned. For these reasons the approximation in a frame is typically more computationally expensive than the approximation in a basis, though for the case of Fourier extensions fast algorithms have been devised \cite{matthysen2016fastfe,matthysen2017fastfe2d}.

\subsection{Computing a frame approximation}

Two ways to compute a frame approximation are investigated in \cite{huybrechs2016tw674}. The first is the best approximation to $f$ by projecting onto the span of the finite set $\Phi_N := \{ \phi_k \}_{k=1}^N$. This leads to the linear system $G^p x = B^p$ where
\begin{equation}\label{eq:continuous}
 G^p_{i,j} = \langle \phi_j, \phi_i \rangle_{L^2(\Omega)} \qquad \mbox{and} \qquad B^p_i = \langle f, \phi_i \rangle_{L^2(\Omega)}.
\end{equation}
The square matrix $G^p \in \mathbb{C}^{N \times N}$ is called the \emph{Gram matrix} in \cite{huybrechs2016tw674}. Both the elements of $G^p$ and of the right hand side $B^p$ require the numerical computation of integrals over the domain $\Omega$, hence we refer to it as a continuous projection method.

The second approach is to compute the best approximation in a discrete, rather than continuous, least squares sense. To that end, one chooses a set $X_M := \{ x_m \}_{m=1}^M$ of $M$ points in $\Omega$. This leads to the linear system $G^d x = B^d$,
\begin{equation}\label{eq:discrete}
 G^d_{i,j} = \phi_j(x_i) \qquad \mbox{and} \qquad B^d_i = f(x_i).
\end{equation}
Compared to the continuous projection, the discrete approach is both simpler and cheaper to implement, as it does not require the computation of integrals. However, the matrix $G^d \in \mathbb{C}^{M \times N}$ is now rectangular. We say that $f$ is \emph{oversampled} if $M > N$.

Both matrices $G^p$ and $G^d$ are ill-conditioned. When a large number of points $M$ is used, the properties of the matrix $G^p$ are similar to those of the normal equations $(G^d)^* G^d$ of the discrete least squares problem. This is common in least squares methods \cite{lawson1996leastsquares}. As a result, the condition number of $G^p$ is almost the square of the condition number of $G^d$.

The systems $G^p c = B^p$ and $G^d c = B^d$ are usually solved with a direct solver. The method investigated in \cite{huybrechs2016tw674} is based on a truncated singular value decomposition (SVD) of the matrix. That is, the systems are regularised by discarding (truncating) all singular values below a threshold $\epsilon$. The error bounds further on are derived for this truncated SVD solution. Numerical experiments show that, typically, slightly higher accuracy is obtained by using a rank-revealing pivoted QR-decomposition.\footnote{We like to stress here that the choice of the direct solver matters. For example, the solution $x = A \setminus B$ using the backslash-operator in MATLAB does not necessarily satisfy the bounds of this paper. The bounds are derived for a solution obtained using a truncated SVD, which yields a solution vector with minimal norm $\Vert x \Vert$. In fact, we found that for several examples in this paper, MATLAB's backslash solution was unstable due to growth of $\Vert x \Vert$ after an initial regime of convergence. In contrast, both the SVD and QR-based solution remain stable for larger systems and hence should be preferred for this class of ill-conditioned matrices.} For this reason, in our implementation, we have used a QR-based solver.

\subsection{Error bounds}

Can one expect highly accurate approximations, in spite of the potentially extreme ill-conditioning of $G$? The surprising answer given in \cite{huybrechs2016tw674} is that, yes, in many cases you can, even when $G$ is nearly singular. An intuitive explanation is that, because of the redundancy of frames, there are many ways in which $f$ can be expressed as a linear combination of the frame elements with comparable accuracy. Thus, it does not matter which solution of $Gx=B$ is found. Instead, the approximation error is governed by the residual $Gx-B$. Achieving a small residual is a considerably weaker condition than finding the exact solution of an ill-conditioned system of equations. Indeed, the ill-conditioning implies that the matrix $G$ has a large numerical kernel, yet any element of the kernel (satisfying $Gx \approx 0$) has no impact on the residual.

Even if $G$ is mathematically non-singular, and thus the system has a unique exact solution, finding this particular solution numerically is not feasible. The \emph{regularized} solution found by truncating the SVD may in fact differ substantially from the exact solution, at least when comparing the coefficients. Yet, in most cases one does not a-priori care which set of coefficients is found, as long as the approximation $\sum_{k=1}^N c_k \phi_k(x)$ is accurate, \emph{and} as long as the coefficients $c_k$ do not grow too large. One typically prefers small coefficients $c_k$, since large coefficients with alternating signs may yield large cancellation errors in subsequent computations.

\emph{These observations lead to two criteria: we are looking for a solution $x$ with small residual $Gx-B$, and with small norm $\Vert x \Vert$.}

A first generic error bound formulated in \cite{huybrechs2016tw674} for the approximation of a function $f$ via the truncated SVD of the projection matrix $G^g$ is:
\begin{equation}\label{eq:error_galerkin}
 \Vert f - P_N^{p,\epsilon} f \Vert_{L^2(\Omega)} \leq \Vert f - \sum_{k=1}^N z_k \phi_k \Vert_{L^2(\Omega)} + \sqrt{\epsilon} \Vert z \Vert, \qquad \forall z \in \mathbb{C}^N.
\end{equation}
Here, the notation $\Vert f - P_N^{p,\epsilon} f \Vert_{L^2(\Omega)}$ in the left hand side represents the approximation error, i.e., the difference between $f$ and its approximation obtained via regularised SVD with truncation parameter $\epsilon$, measured in the $L^2$ norm on $\Omega$. The inequality \eqref{eq:error_galerkin} holds for any vector $z$, and hence shows that the solution will be at least as good as any approximation to $f$ in the span of the truncated frame $\Phi_N$, as long as the size of the coefficients $\Vert z \Vert$ remains small. One finds a best approximation, subject to having small coefficients. The maximal achievable accuracy is on the order of $\sqrt{\epsilon}$.

A second generic error bound applies to the oversampled least squares problem \eqref{eq:discrete}.  It is slightly more involved and introduces more notation:
\begin{align}\label{eq:error_collocation}
 \Vert f - P_{M,N}^{d,\epsilon} f \Vert_{L^2(\Omega)} \leq \Vert f - \sum_{k=1}^N z_k \phi_k \Vert_{L^2(\Omega)} &+ \kappa_{M,N}^\epsilon \, \Vert G^d z - B^d \Vert \nonumber \\ 
 &+ \epsilon \, \lambda_{M,N}^{\epsilon} \, \Vert z \Vert, \qquad \forall z \in \mathbb{C}^N.
\end{align}
As above, the quantity in the left hand side is the approximation error, measured on $\Omega$. The right hand side also features the residual of the linear system corresponding to vector $z$. 
In addition, there are two constants, $\kappa_{M,N}^{\epsilon}$ and $\lambda_{M,N}^{\epsilon}$. These are difficult to analyse generically. However, in the case where $\Phi$ is a frame, these constants are bounded for sufficiently large $M$ (with $M$ large compared to $N$). This corresponds to sufficiently high oversampling.

Crucially, the maximal achievable error of the oversampled discrete scheme is on the order of $\epsilon$ in \eqref{eq:error_collocation}, rather than $\sqrt{\epsilon}$ in \eqref{eq:error_galerkin}. This is the main reason for the observation that, when manipulating expansions that are possible redundant, it is worthwhile to introduce oversampling into the discretisation of the problem.

Results are also included in \cite{huybrechs2016tw674} on the size of the solution vector, which we do not repeat here. Briefly, before the onset of convergence, the norm of the solution vector may in fact be as large as $1/\epsilon$ or $1/\sqrt{\epsilon}$. However, once $N$ is sufficiently large such that the approximation starts to converge, the norm of the expansion coefficients decreases to a value that can be related to $\Vert f \Vert_{L^2(\Omega)}$, the continuous norm of $f$ itself on the domain $\Omega$. Thus, after potential initial divergence, the regularised computation leads to a numerical approximation that is well-behaved for postprocessing purposes, i.e., it does not suffer from large cancellation errors. In view of the ill-conditioned linear system, this is a very desirable outcome.

\section{A collocation approach of WBM}
\label{sect:collocation}

In the restricted setting of this paper, i.e., in the absence of a multilevel approach or hybrid couplings, the Wave Based Method can be seen as an approximation problem: one approximates a function on the boundary $\partial \Omega$ of a domain $\Omega$, using a set of functions defined on a larger bounding box. Indeed, compare the linear system of the weighted residual formulation \eqref{eq:integrals} to the Gram matrix in \eqref{eq:continuous}. If the expansion succeeds in matching the given Neumann data then, by construction, the expansion also satisfies the Helmholtz problem and a solution to the boundary value problem \eqref{eq:Helmholtz}--\eqref{eq:boundary} is found.

Two drawbacks of the WBM are that (i) it relies on accurate integration techniques and (ii) it seems restricted to convex domains. The collocation formulation of the WBM starts by addressing the former one. The idea is straightforward, by extending the comparison to the frame approximation problem and considering a discrete least squares variant. To that end, the boundary condition is enforced in a set of $M$ points distributed along the boundary, referred to as collocation points $\{s_k\}_{k=1}^M$. This leads to a linear system of equations $\tilde{A}\alpha=\tilde{b},$ where $\tilde{A}$ is an $M\times N$ matrix and $b$ is a vector, with entries given by
\begin{equation}
\label{eq:collocation_system}
\tilde{a}_{i,j} = \phi_j(s_i) \quad \mbox{and} \quad  \tilde{b_i} = v(s_i), \qquad i = 1,\ldots, M,\,\,j = 1,\ldots, N.
\end{equation}
These are the expressions for a Dirichlet problem. The analogous formulation for Neumann boundary conditions is similar, but using the normal derivatives of the functions $\phi_j$ instead.

Bearing in mind the similarities to the discrete least squares approximation problem \eqref{eq:discrete}, it is important here to oversample. That is, we use more collocation points than there are degrees of freedom by a factor $\gamma$, $M=\gamma N$ with $\gamma>1$. The theory of frames does not generically stipulate what $\gamma$ should be, nor that such linear oversampling is even sufficient.\footnote{By linear oversampling we mean the linear relation $M= \gamma N$ between $M$ and $N$. Indeed there are known cases where $M$ should grow like $N^2$. Needless to say, the amount of oversampling negatively effects the computational cost of the scheme.} The value of $\gamma$ is problem dependent. In this paper, unless stated otherwise, we shall use a value of $\gamma = 2$ and we observe numerically that this appears sufficient for our experiments.

A first immediate benefit from this approach is that the collocation matrix $\tilde{A}$ is less severely affected by ill-conditioning than the WBM matrix $A$. This follows from similar considerations discussed for frame approximations. A second benefit is computational: the cost of assembling $\tilde{A}$ is significantly lower, since each element of the matrix consists of one function evaluation, whereas in the weighted residual formulation each element requires the numerical approximation of an integral over the boundary of the domain.

While our formulation of the oversampled collocation formulation is based on the connection to frames, we like to point out that other authors have made similar implementations for Trefftz methods. For example, the method of fundamental solutions as implemented in \cite{barnett2008stability} is based on oversampled collocation.

\section{Results}
\label{sect:results}

We apply both the weighted residual and the oversampled collocation formulations of WBM for a number of two-dimensional Helmholtz problems on smooth domains. We use either Dirichlet or Neumann boundary conditions. We focus on two properties: convergence of the scheme as a function of the number of degrees of freedom, and size $\Vert x \Vert$ of the coefficient vector in the discrete $l^2$-norm.

To that end, we devise a number of scattering configurations for which we can show a priori whether or not a solution vector with bounded norm exists. Inspired by the analysis in \cite{barnett2008stability}, the problems involve two kinds of singularities: singularities induced by the boundary data, and singularities induced by the (non-convex) shape of the domain. For the time being, we do not consider singularities due to corners or other singularities of the domain shape itself.

All experiments were performed in MATLAB. The linear systems for the weighted residual formulation are solved with a custom rank-revealing QR-decomposition with column pivoting, with threshold manually set to $2 \times 10^{-13}$, a value that was experimentally determined. The collocation systems were solved simply using MATLAB's backslash operator, which runs a comparable algorithm for this class of matrices. The integrals \eqref{eq:integrals} involved with the classical WBM system are approximated using the trapezoidal rule with sufficiently many points, since all integrands are smooth and periodic \cite{trefethen2014trapezoidal}.

Since the wave functions satisfy exactly the Helmholtz equation inside the domain, errors are only introduced by not exactly satisfying the boundary conditions. Though these errors propagate towards the interior of the domain, we measure the relative error on the boundary data in a set of points $\{ (x_i,y_i) \}_{i=1}^{n_p}$ on $\partial \Omega$, with $n_p$ larger than the number of collocation points used to compute the approximation. For example, for Dirichlet data, this results in the metric
\begin{equation}
\label{eq:error}
\varepsilon = \frac{1}{n_p}\left(\sum\limits_{i = 1}^{n_p}\left|\frac{\tilde{u}(x_i,y_i)-w(x_i,y_i)}{w(x_i,y_i)} \right| \right),
\end{equation}
where $w$ denotes the function imposed on the boundary. 

In several experiments, we choose boundary conditions that correspond to a known solution of the Helmholtz equation that is bounded in the interior of $\Omega$. Hence, by construction, we know the exact solution of the interior boundary value problem and, moreover, we know the extension of this solution outside the domain $\Omega$.

In this section we describe the numerical experiments and the results. We elaborate on their interpretation further on in \S\ref{sect:discussion}.

\subsection{Convex obstacles}

We consider convex obstacles first and start by examining the effect of the size of the bounding box. This experiment is inspired by the analogy to Fourier extension frames, illustrated in Fig.~\ref{fig:extensionframe} which should be compared with Fig.~\ref{fig:WBMdomains}. For Fourier extension approximations, faster convergence rates are seen for larger bounding boxes \cite{huybrechs2010fourier,adcock2014resolution}.

\begin{figure}[tbp]
 	\centering
 	\begin{subfigure}[c]{0.47\textwidth}
 		\includegraphics[width=\textwidth]{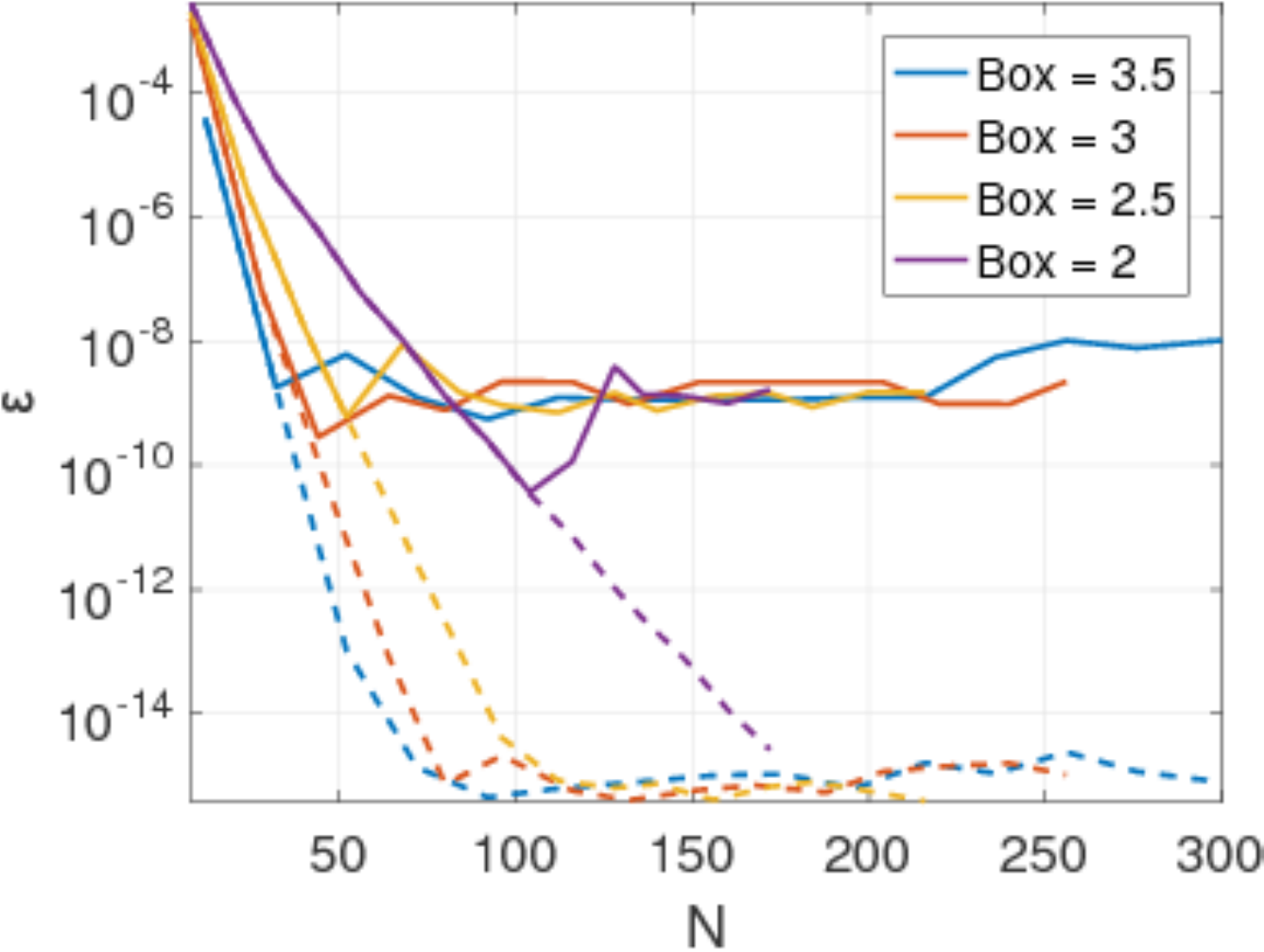}
 		\caption{Convergence}
 		\label{fig:BoxConvg}
 	\end{subfigure}
 	\begin{subfigure}[d]{0.47\textwidth}
 		\includegraphics[width=\textwidth]{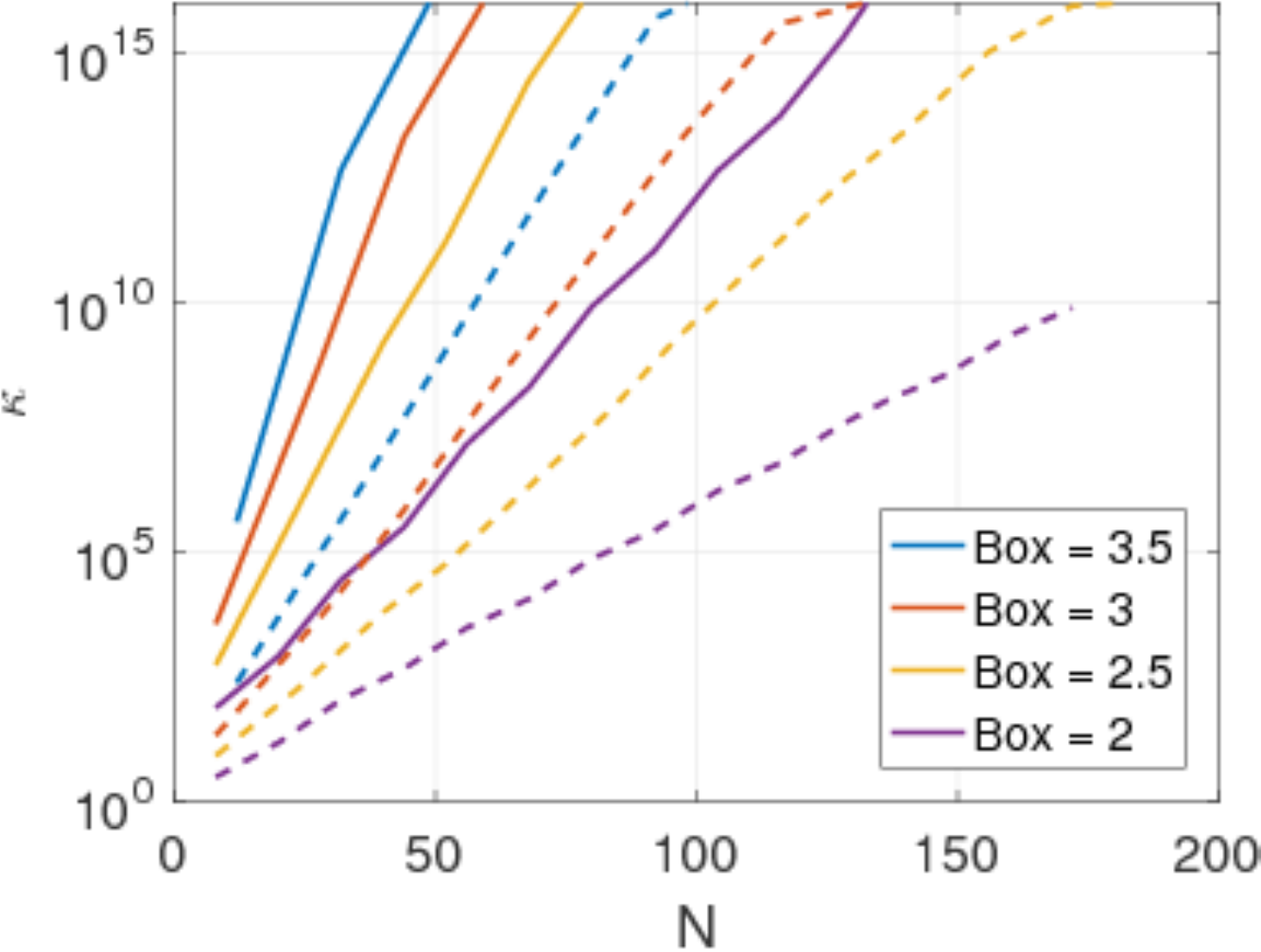}
 		\caption{Condition number}
 		\label{BoxCond}
 	\end{subfigure}   	        
 	\begin{subfigure}[d]{0.47\textwidth}
 		\includegraphics[width=\textwidth]{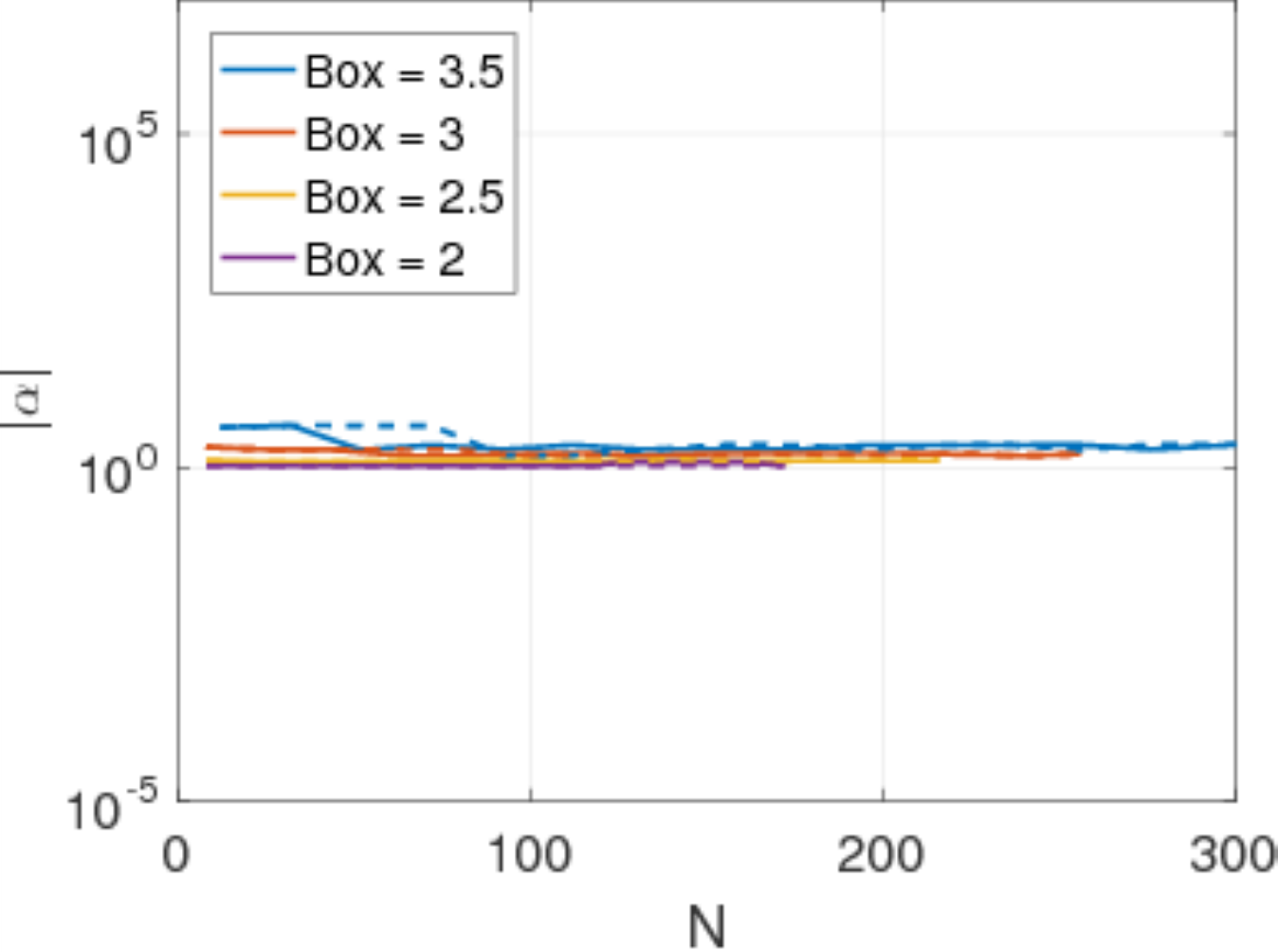}
 		\caption{Coefficient norm}
 		\label{BoxCoeff}
 	\end{subfigure}   		
 	\caption{The solution on the unit disk with wavenumber $k=0.927$, embedded in a square bounding box with edge lengths varying between $2$ and $3.5$. Solid lines correspond to weighted residual WBM, dashed line to the collocation approach.}
 	\label{fig:WaveBox}
 \end{figure}

The domain is the unit disk centred at the midpoint of a square bounding box. We choose a wavenumber $k=0.924$ and Neumann boundary conditions consisting of the normal derivative of a plane wave propagating at an angle of $0.3$ radians. Thus, the exact solution on the disk is this plane wave. Results are shown in Fig.~\ref{fig:WaveBox} as a function of $N$, the total number of degrees of freedom in the discretisation, and for different sizes of the edge length of the bounding square. The results show the accuracy, the condition number of the linear system and the $l^2$-norm of the coefficient vector.

One observes in Fig.~\ref{fig:WaveBox}(a) that the collocation approach achieves higher accuracy than the weighted residual formulation, on the order of $10^{-15}$ compared to $10^{-8}$. Higher convergence rates are seen for larger box sizes. Panel (b) shows very large condition numbers, and confirms that $\kappa(A) \approx \kappa(\tilde{A})^2$. Finally, panel (c) shows that the computed solution vectors have a norm that remains bounded as $N$ increases, for both variations of WBM.



Next, we consider the same circular domain centred at $(1.5,1.5)$, and enclosed in a $[0,3]\times[0,3]$ box of fixed size. This time we employ Dirichlet boundary conditions arising from a function satisfying the Helmholtz equation in the plane, yet with a singularity at a point $(x_s,y_s)$. In particular, we use the Hankel function of the first kind and order zero,  $w(x,y) = H_0^{(1)}(k\Vert (x,y)-(x_s,y_s)\Vert)$. We examine convergence and norm of the solution when varying the location of the singularity, which corresponds to a point source. The ordinate $y_s = 1.5$ is kept constant, while $x_s \in [-1,\, -0.01,\, 0.01,\, 0.2,\, 0.4 ]$.

\begin{figure}[tbp]
  	\centering
  	\begin{subfigure}[c]{0.45\textwidth}
  		\includegraphics[width=\textwidth]{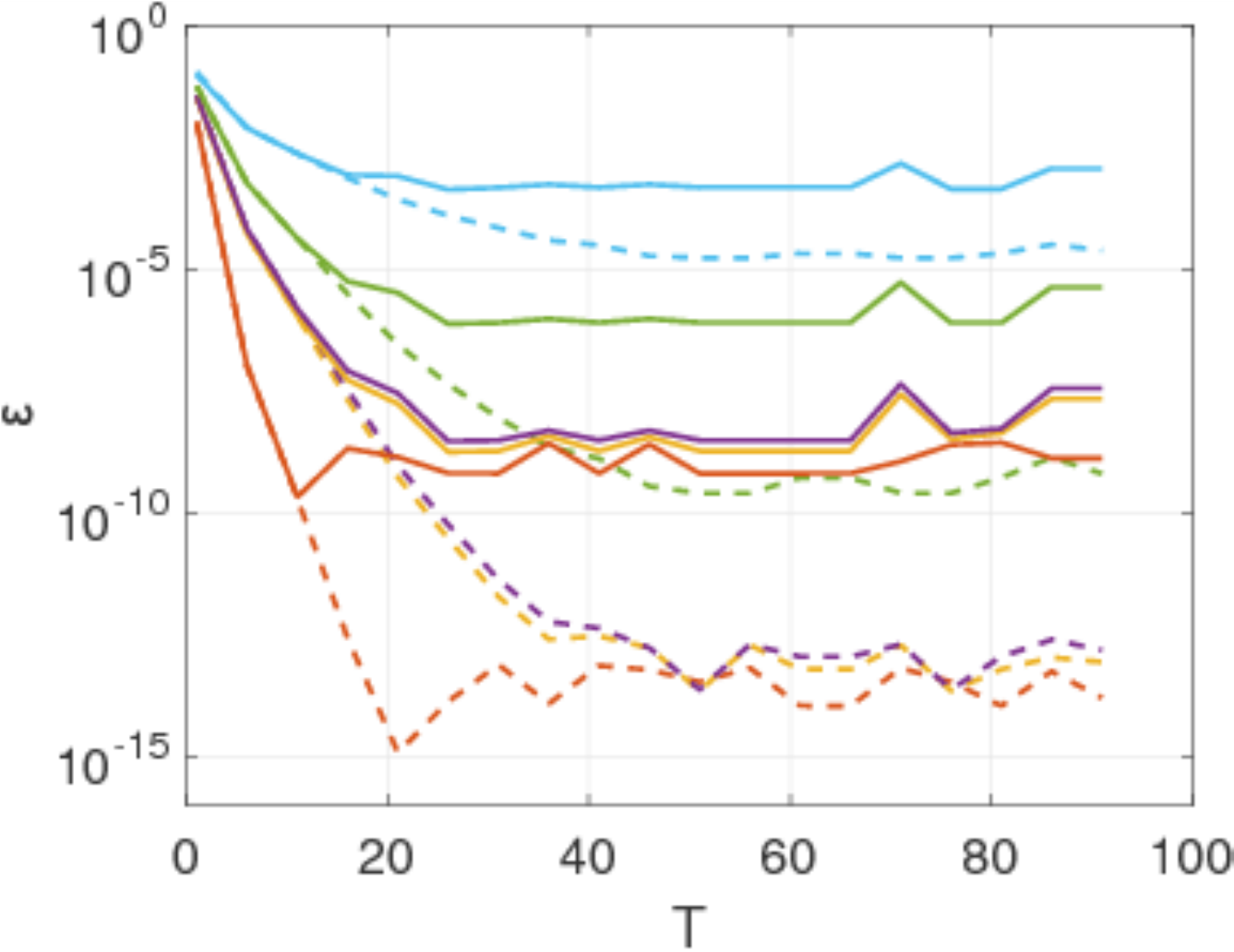}
  		\caption{Convergence}
  		\label{fig:CircleHankelErr}
  	\end{subfigure}     	
  	\begin{subfigure}[d]{0.45\textwidth}
  		\includegraphics[width=\textwidth]{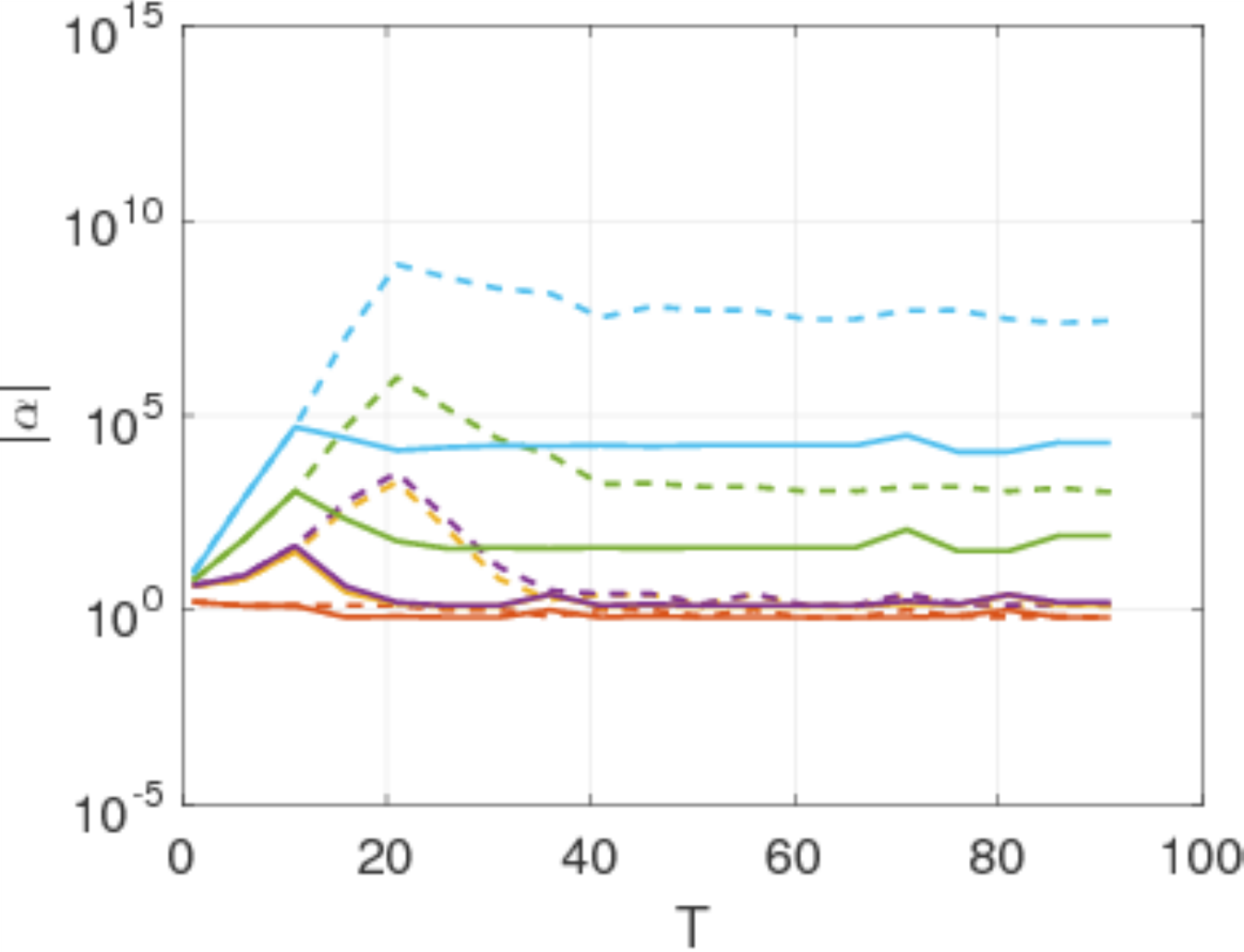}
  		\caption{Coefficient norm}
  		\label{fig:CircleHankelCoeff}
  	\end{subfigure}       	
  	\caption{WBM results using the weighted residual (line) and collocation (dashed line) formulations on a disk with point source boundary data; color varies with the abscissa of the source point: $-1$ (red), $-0.01$ (yellow), $0.01$ (purple), $0.2$ (green), $0.4$ (light blue).}
  	\label{fig:CircleHankel}
  \end{figure}

The results are shown in Fig.~\ref{fig:CircleHankel}. Here, too, the collocation approach converges to higher accuracy as shown in panel (a). However, both implementations appear to have a maximal accuracy that depends on the location of the singularity. When the singularity is furthest away from the box ($x_s = -1$), the convergence curves are comparable to those in Fig.~\ref{fig:WaveBox}(a). As the source approaches the box, both outside ($x_s = -0.01$) and inside ($x_s = 0.01$) the box, a smaller convergence rate is seen, yet similar accuracy is ultimately achieved. However, in panel (b), a small bump is seen in the corresponding curves for small $T$: the solution vector $\Vert x \Vert$ initially grows, before settling down for increasing $T$ to a small value comparable to the first case of a far away singularity. The bump is larger for the collocation approach.

The impact of the singularity worsens as it approaches the disk even more. For $x_s=0.2$ and $x_s=0.4$, the convergence curves in Fig.~\ref{fig:CircleHankel} level off at larger and larger values. It remains the case that collocation reaches higher accuracy than the weighted residual method. However, this comes at a price, seen in panel (b): the norms of the solution vectors settle down at significantly larger levels. In addition, the norms are significantly larger in the case of collocation.

As the singularity approaches the domain, both methods become increasingly inaccurate and increasingly unstable. The best accuracy achieved for both formulations of the WBM is restricted by the growth of the norm of the coefficients. In this case, convergence stops when $\frac{\varepsilon}{\alpha}\approx 10^{-13}$ for collocation, and $\frac{\varepsilon}{\alpha}\approx 10^{-8}$ for the weighted residual approach, where $\alpha$ is the limiting value of $\Vert x \Vert$ for large $T$. These regimes corresponds precisely to an equal balance of the terms appearing in the right hand sides of \eqref{eq:error_galerkin} and \eqref{eq:error_collocation}, respectively.

\subsection{Non-convex obstacles}

In the previous experiment, the solution of the boundary value problem has a singularity in the exterior of the domain $\Omega$, determined by the boundary condition. It is known that the solution of a Helmholtz problem on a non-convex domain also develops singularities in its extension to the exterior. However, in this case, the locations of the singularity are determined by the geometry of the domain, in particular its boundary, and not by the boundary condition data. We refer the reader to \cite{barnett2008stability} and the references therein for an extensive description of known results. We adopt some examples from the same reference to illustrate the impact of these singularities in the Wave Based Method.

\begin{figure}[tbp]
	\centering
	\begin{subfigure}[c]{0.32\textwidth}
		\includegraphics[width=\textwidth]{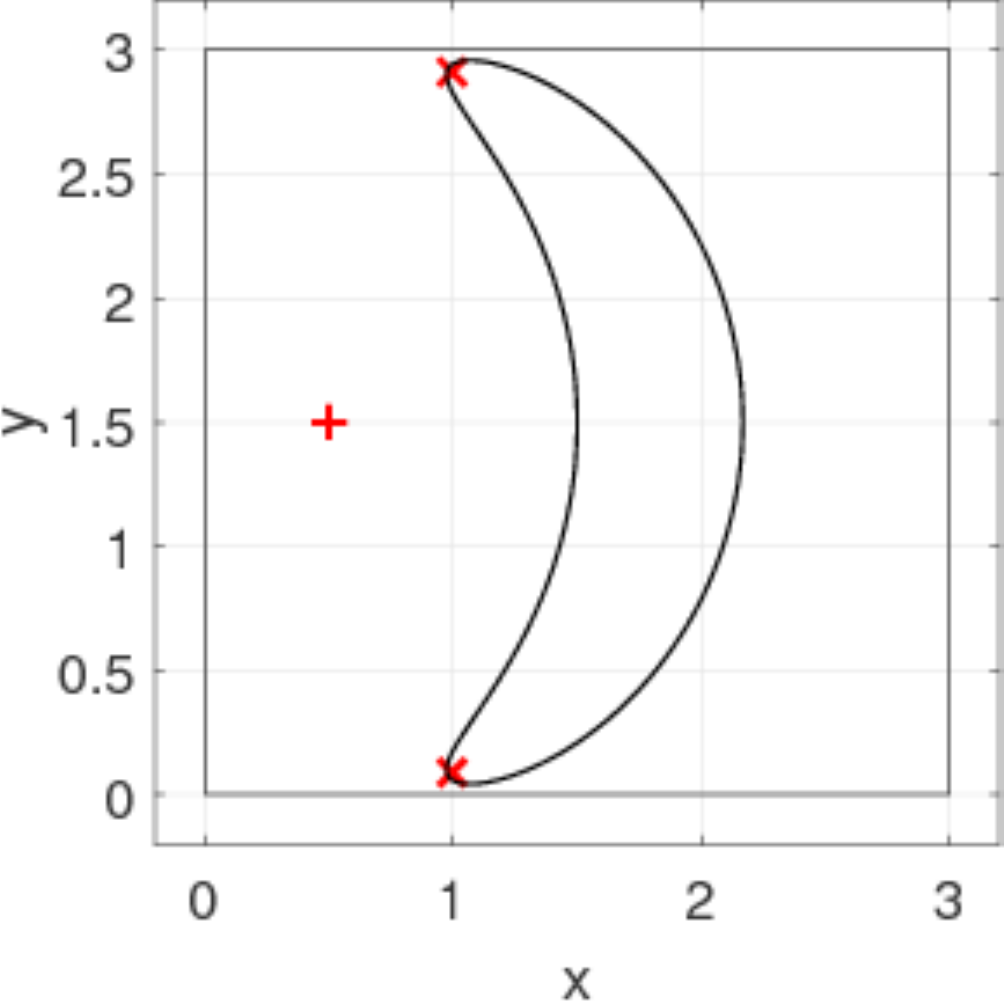}
		\caption{ \Romannum{1}: $a = b = 0.5$}
		\label{fig:cresDomI}
	\end{subfigure}     	
	\begin{subfigure}[d]{0.32\textwidth}
		\includegraphics[width=\textwidth]{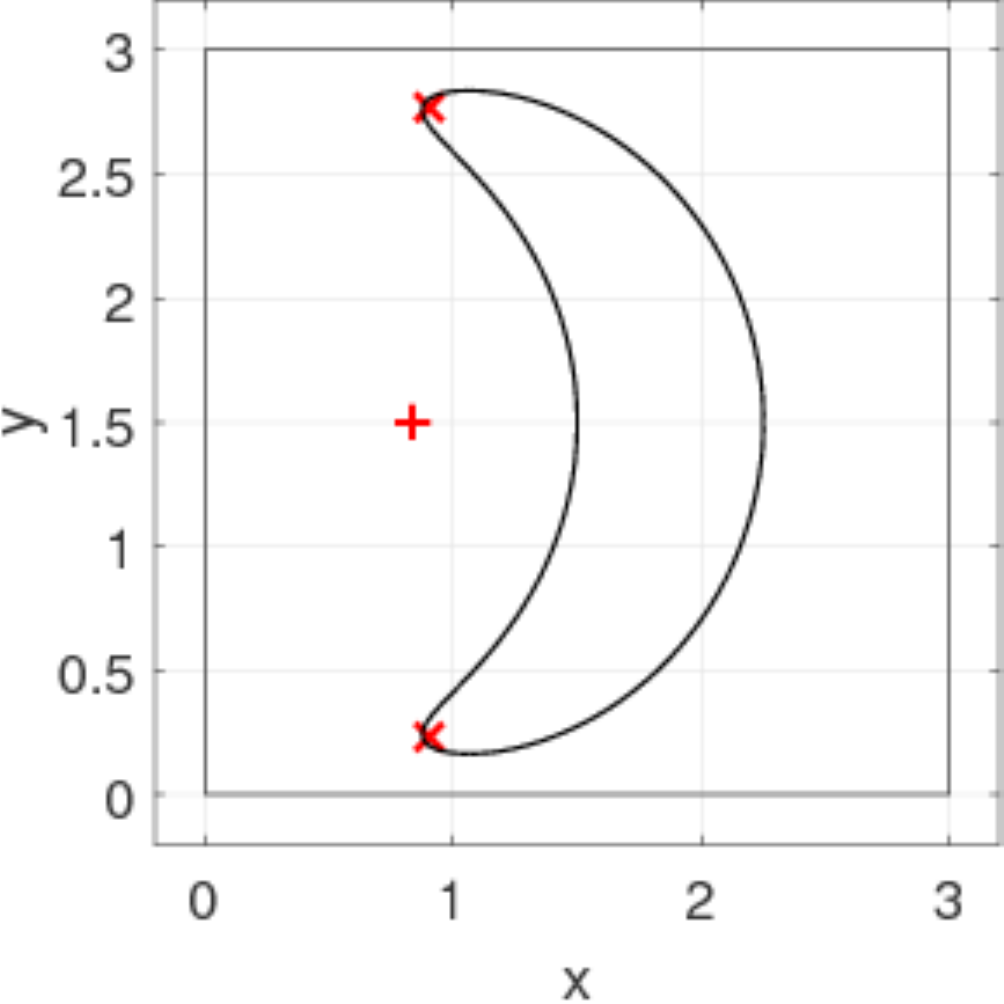}
		\caption{ \Romannum{2}: $a = 0.4,\,\,b=0.6$}
		\label{fig:cresDomII}
	\end{subfigure}   
	\begin{subfigure}[d]{0.32\textwidth}
		\includegraphics[width=\textwidth]{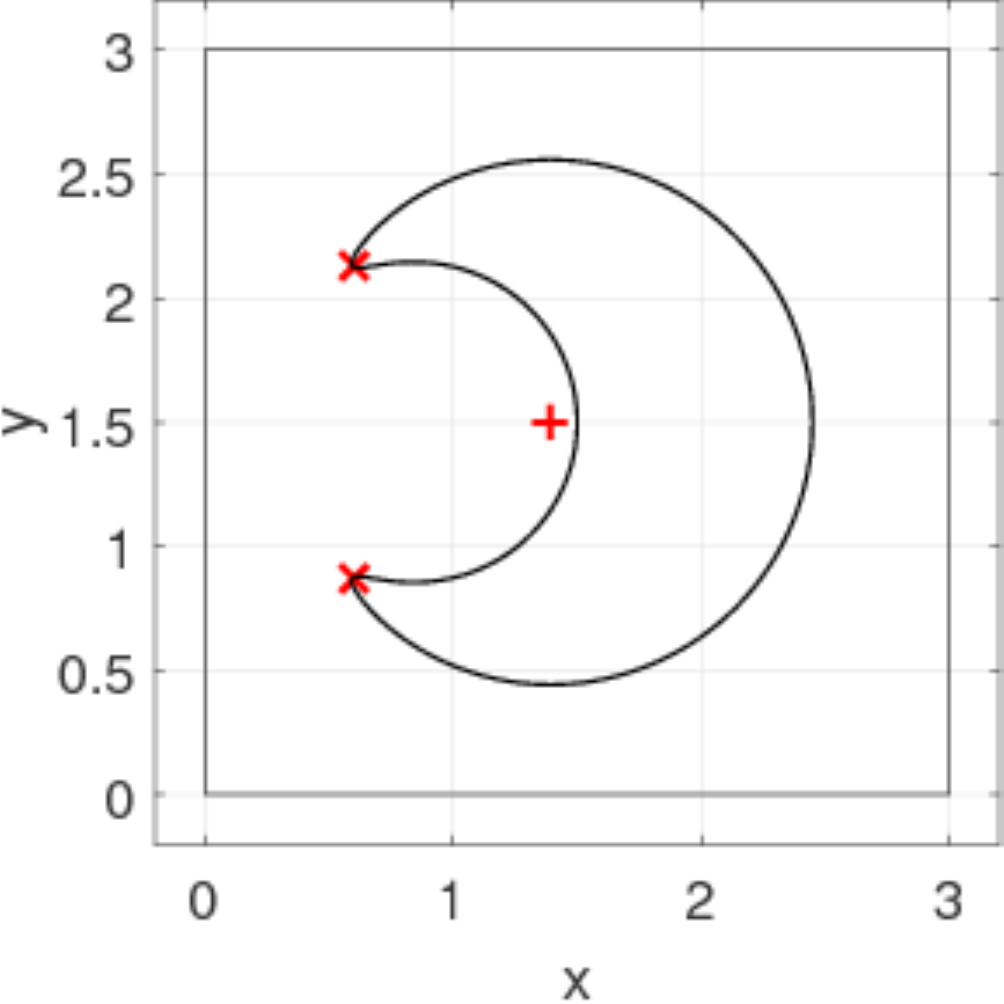}
		\caption{ \Romannum{3}: $a = 0.1, b = 0.9$}
		\label{fig:cresDomIII}
	\end{subfigure}      	
	\caption{The relevant singularities for a non-convex domain are those of the associated Schwartz function. Shown here is the crescent domain of \cite{barnett2008stability} and its singularities: a pole in the exterior, denoted by '$+$', and two branch point singularities in the interior denoted by '$\times$. The parameterisation of the boundary is given by \eqref{eq:cresParam} and contains two parameters $a$ and $b$. In all cases shown, the crescent domain is enclosed by a square with edge length $3$.}
	\label{fig:cresDom}
\end{figure}

We show in Fig.~\ref{fig:cresDom} the domain referred to as `crescent' in \cite{barnett2008stability}. The parameterisation of its boundary is given by
 \begin{equation}\label{eq:cresParam}
f_1(t) = z_0 + e^{it}- \displaystyle\frac{a}{e^{it}+b}, \qquad t \in [0,2\pi],
 \end{equation}
where $z_0 = 1.5+1.5i$ and $a$ and $b$ are two parameters. Here, $f_1(t) = x(t) + i y(t)$ is a complex-valued function and we identify its real and imaginary parts with coordinates in the Euclidean plane $\mathbb{R}^2$.

The singularities associated with a non-convex domain that has an analytic parameterisation (analytic in the sense of complex analysis) are those of the so-called Schwartz function. For a given domain $\Omega$ with boundary $f(t)$, where $f$ is an analytic function of $t$ and $|f'(t)| \neq 0$, the Schwartz function is defined as $S(z) =  \overline{f}(f^{-1}(z))$. It is a function of $z$ in the complex plane, and may have singularities such as branch points and poles. We omit their computation, but highlight for the crescent domain in Fig.~\ref{fig:cresDom}(a)--(c) the location of a pole in the exterior of $\Omega$, and of two branch points of square root type in the interior of $\Omega$, for different combinations of the parameters $a$ and $b$. The branch points are located at $z_{b_{1,2}} =z_0-b\pm2i\sqrt{a}$ and the pole at $z_p = z_0-\displaystyle\frac{a}{b}$.

\begin{figure}[!htbp]
 	\centering
 	\begin{subfigure}[c]{0.45\textwidth}
 		\includegraphics[width=\textwidth]{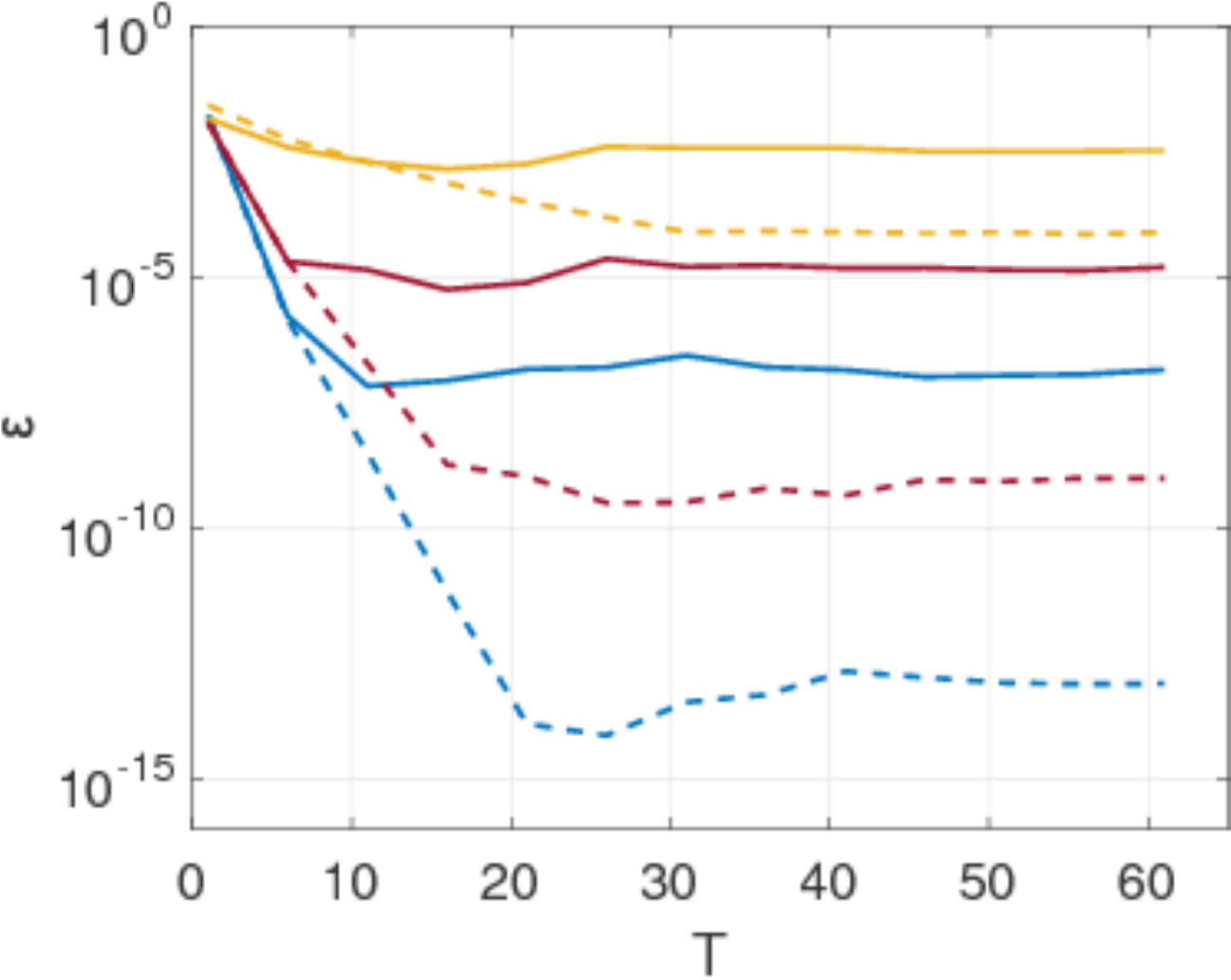}
 		\caption{Approximation error}
 		\label{fig:CresConstConvg}
 	\end{subfigure}
 	\begin{subfigure}[d]{0.45\textwidth}
 		\includegraphics[width=\textwidth]{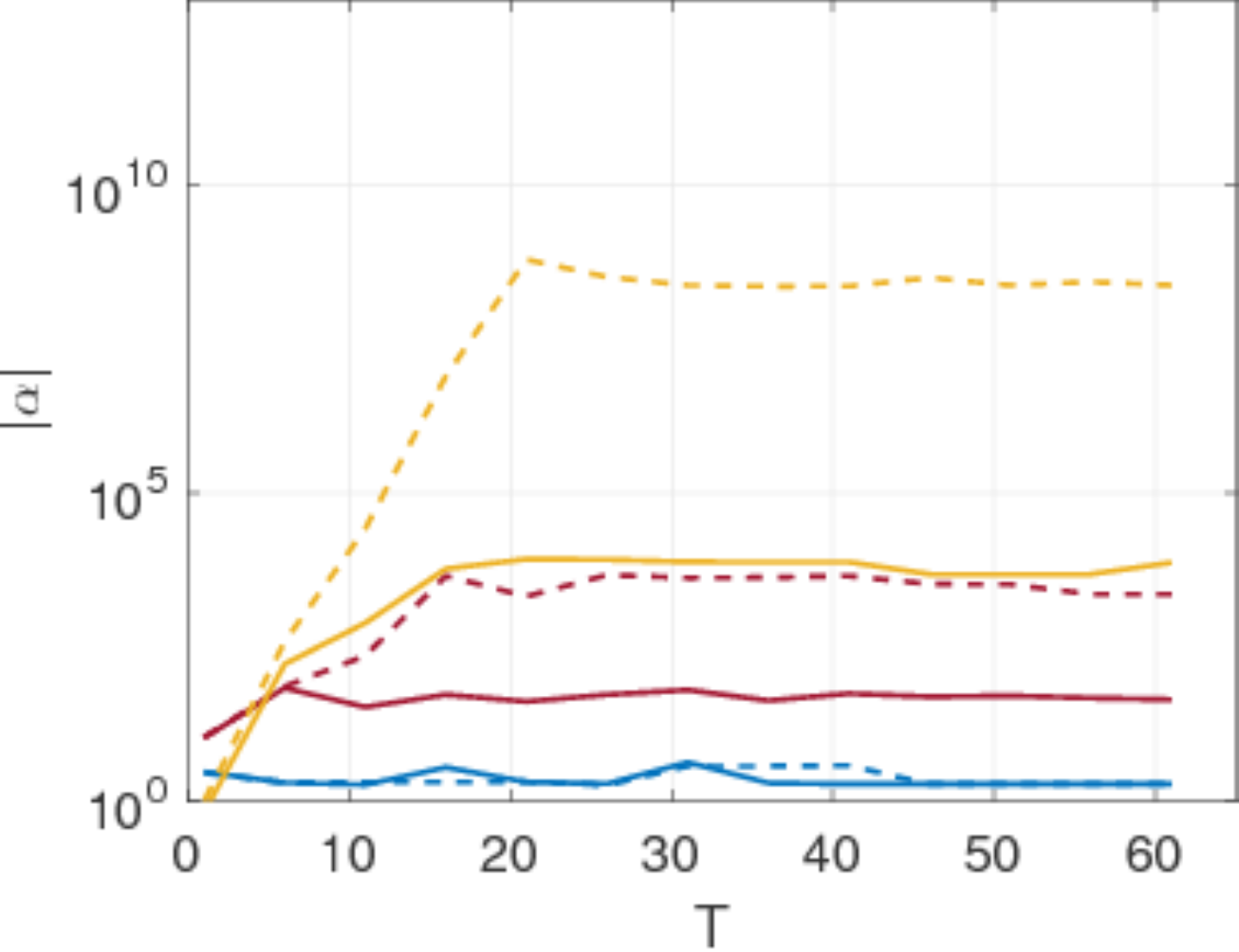}
 		\caption{$2-$norm of the coefficient vector}
 		\label{fig:CresConstCoeff}
 	\end{subfigure}   	   	
 	\caption{WBM results using the weighted residual (line) and collocation (dashed line) formulations with constant Dirichlet boundary data $w(x,y)=1$, for the three variations of the crescent domain shown in Fig.~\ref{fig:cresDom}: \Romannum{1} (blue), \Romannum{2} (red), \Romannum{3} (yellow). } 
 	\label{fig:CresConst}
 \end{figure}

The solution of equation $\ref{eq:Helmholtz}$ with constant Dirichlet boundary condition, $w(x,y) = 1,$ for $(x,y)\in\partial\Omega$, is approximated for each variation of the crescent domain. The results are shown in Fig.~\ref{fig:CresConst}. The error is calculated using $n_p=3000$ boundary points.

As the domain becomes more non-convex, which corresponds to the pole of the Schwartz function moving closer to the domain, the maximal accuracy of WBM decreases. Similar to the case of convex domains with nearby singularities, the collocation approach yields higher accuracy, at the cost of producing a solution with larger coefficient norm.

\begin{figure}[!htbp]
	\centering
	\begin{subfigure}[c]{0.45\textwidth}
		\includegraphics[width=\textwidth]{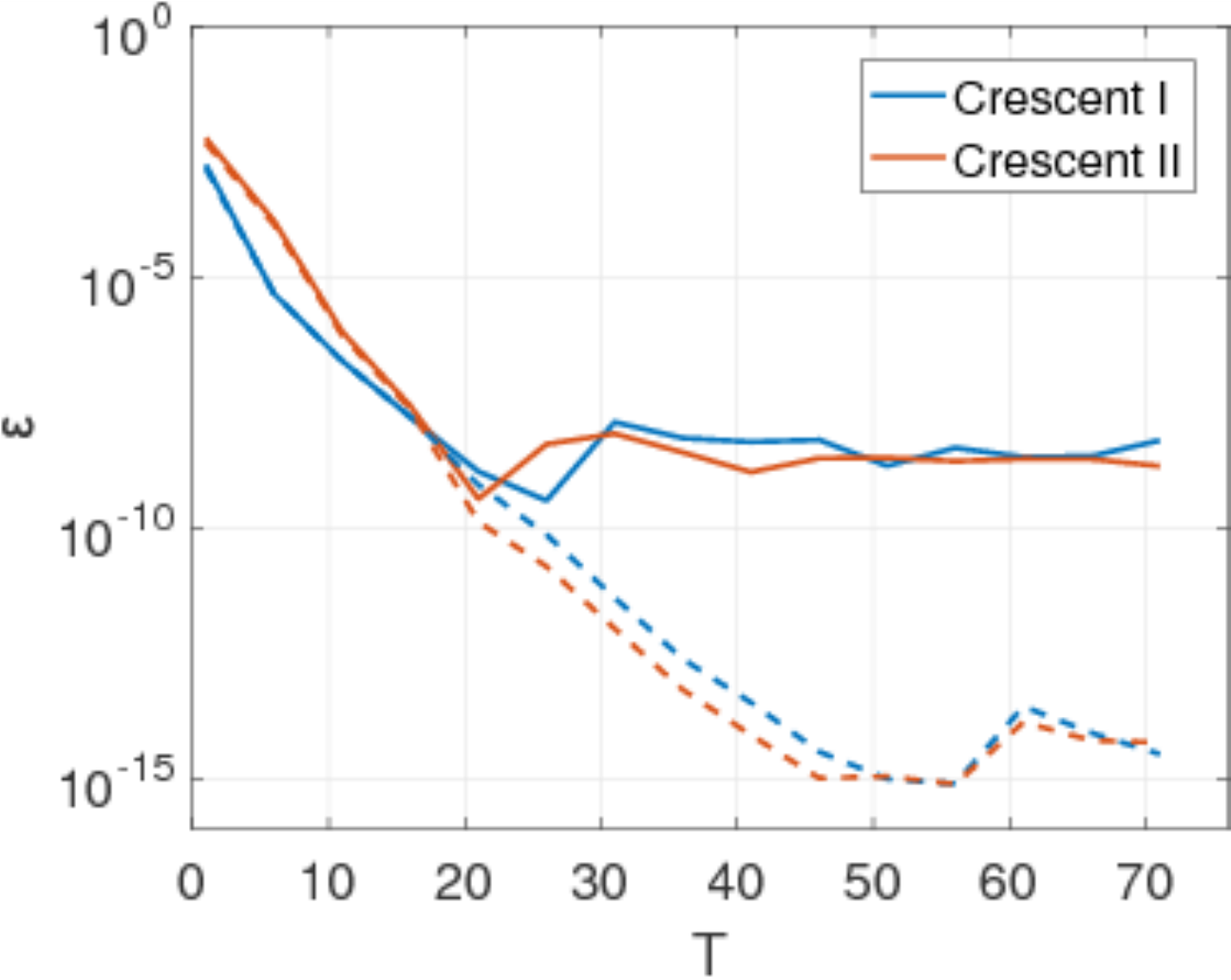}
		\caption{Approximation error}
		\label{fig:CresConstConvgII}
	\end{subfigure}
	\begin{subfigure}[d]{0.45\textwidth}
		\includegraphics[width=\textwidth]{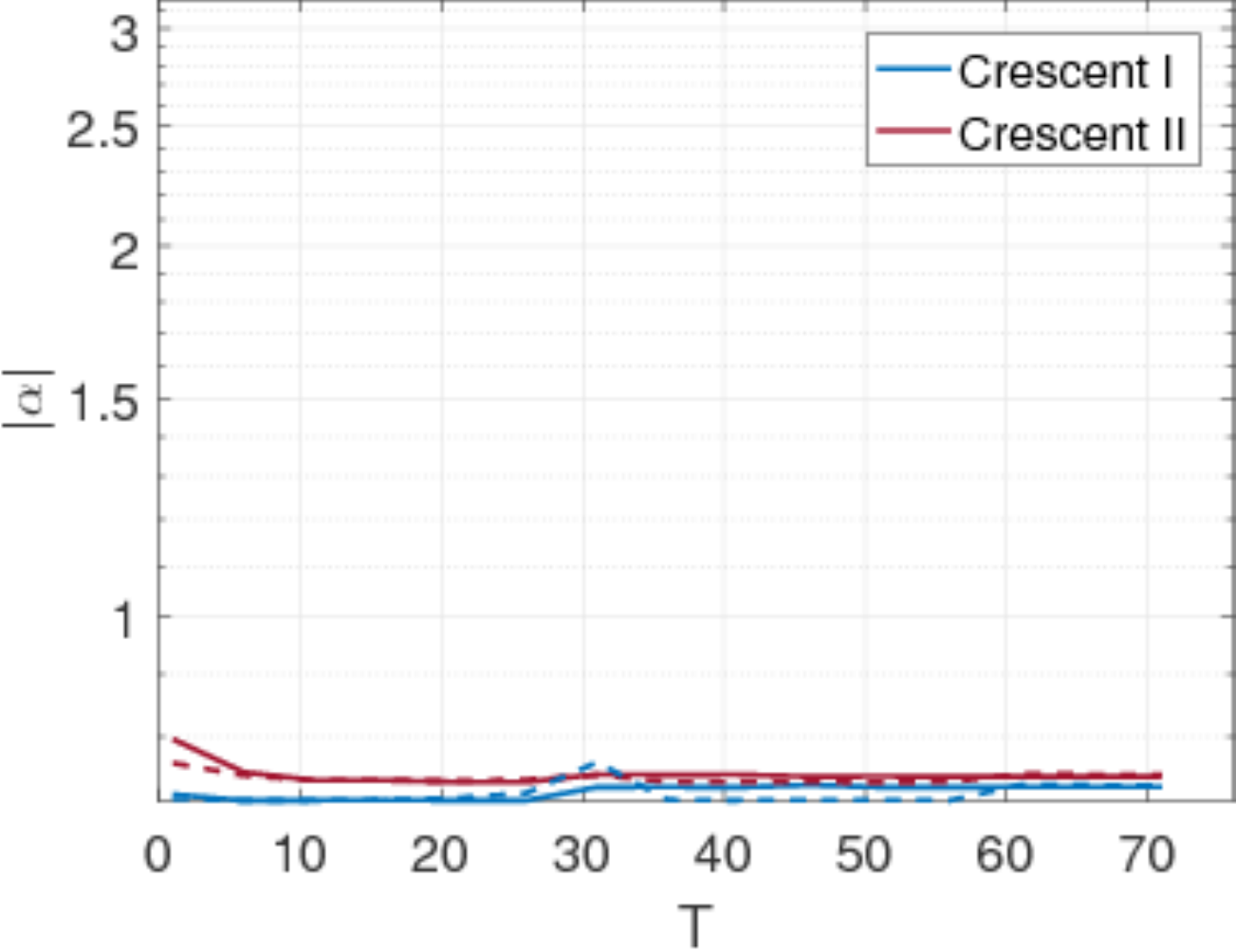}
		\caption{$2-$norm of the coefficient vector}
		\label{fig:CresConstCoeffII}
	\end{subfigure}   	   	
	\caption{WBM results using weighted residual (line) and collocation (dashed line) formulations with constant Dirichlet boundary data $w(x,y)=1$, for crescents \Romannum{1} (blue) and \Romannum{2} (red), and with a bounding box chosen such that the singularities are located outside.}
	\label{fig:CresConstII}
\end{figure}

In the case of crescents \Romannum{1} and \Romannum{2}, it is possible to choose a bounding box that encloses the domain but not the poles of the Schwartz function. In the next experiment the box is centred at $z_0= 0.64+1.75i$, with sides of length $L_x = 1.4$ and $L_y = 3.5$. The results in Fig.~\ref{fig:CresConstII} show that in both cases full expected accuracy is achieved (of order $\sqrt{\epsilon}$ and order $\epsilon$, respectively, for weighted residual and collocation formulations). This is in spite of the non-convexity of the domain and in spite of the ill-conditioning of the linear systems. The coefficient vector has small norm for a wide range of the truncation parameter $T$. The convergence rate is lower compared to the results in Fig.~\ref{fig:CresConst}, because the bounding box is tighter. This is in agreement with the first experiment of this section, in which we varied the size of the box.

This experiment suggests that the convexity requirements in the Wave Based Method is not a fundamental restriction. Rather, at least in the current setting of smooth obstacles, convergence and stability are simultaneously hampered by the presence of singularities in the extension of the Helmholtz solution, exterior to $\Omega$ but within the bounding box. We elaborate on this point in \S\ref{sect:discussion}.

\begin{figure}[!htbp]
 	\centering
 	\begin{subfigure}[c]{0.3\textwidth}
 		\includegraphics[width=\textwidth]{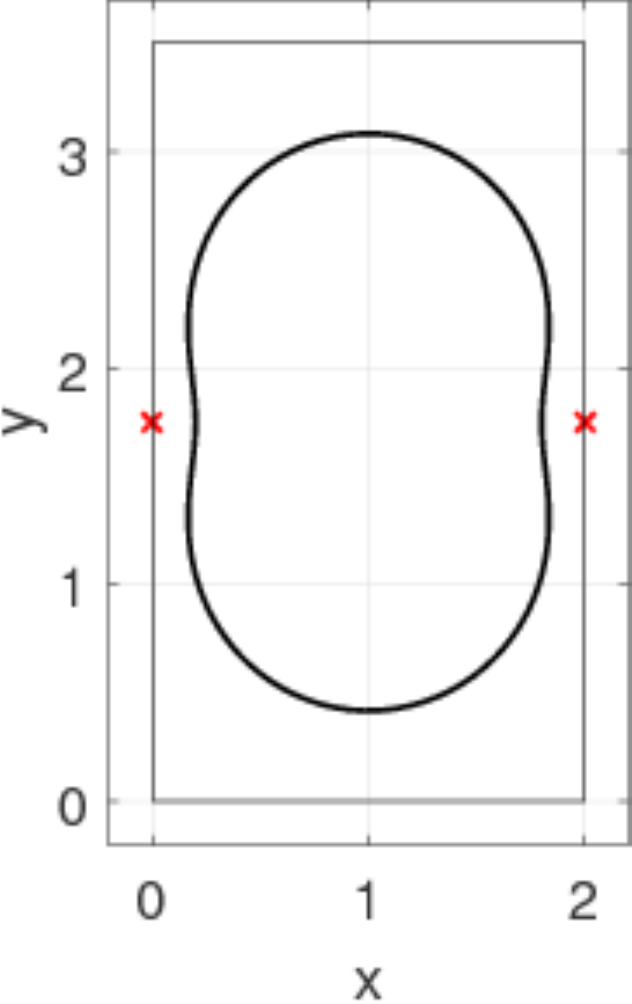}
 		\caption{$\tau = 0.25$}
 		\label{fig:InvElI}
 	\end{subfigure}
 	\hspace{1cm}
 	\begin{subfigure}[d]{0.3\textwidth}
 		\includegraphics[width=\textwidth]{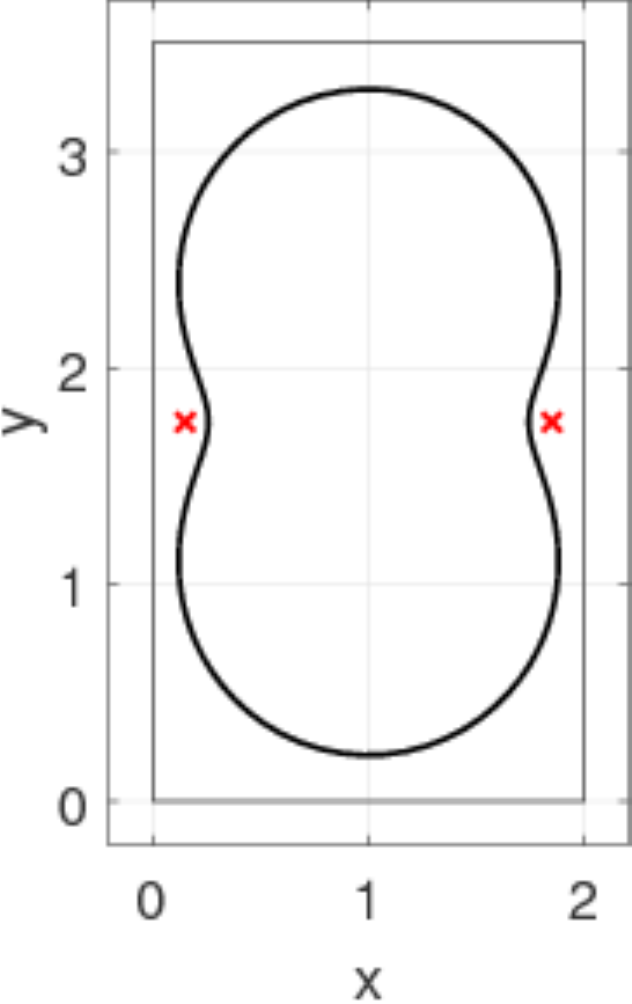}
 		\caption{$\tau = 0.35$}
 		\label{fig:InvElII}
 	\end{subfigure}    	   	
 	\caption{Two inverted ellipses and the location of the two branch points of the corresponding Schwartz functions: they are on the box boundary (left panel) or in its interior (right panel).}
 	\label{fig:InvElDomain}
 \end{figure}

We further illustrate this property with another example domain from \cite{barnett2008stability}, the `inverted ellipse'. The domain is illustrated in Fig.~\ref{fig:InvElDomain}. Its boundary is parameterised by the complex-valued function
\[
f_2(t) = z_0+\frac{e^{it}}{1+\tau e^{2it}}, \qquad t \in [0,2\pi],
\]
where we choose $z_0 = 1 + 1.75i$ and $\tau$ is a parameter. 

We consider two values $\tau = 0.25$ and $\tau = 0.35$ and choose the bounding box to be $[0,2]\times[0,3.5]$. Computation of the Schwarz function shows two branch point singularities at $z_{b_{1,2}}= z_0\pm\sqrt{\frac{1}{4\tau}}$. For the first ellipse ($\tau=0.25$), the singularities are exactly on the edge of the box, whereas in the second case ($\tau=0.35$) they are located in the interior of the box. Indeed, the second inverted ellipse deviates more from a convex domain than the first one.

\begin{figure}[!htbp]
  	\centering
  	\begin{subfigure}[c]{0.45\textwidth}
  		\includegraphics[width=\textwidth]{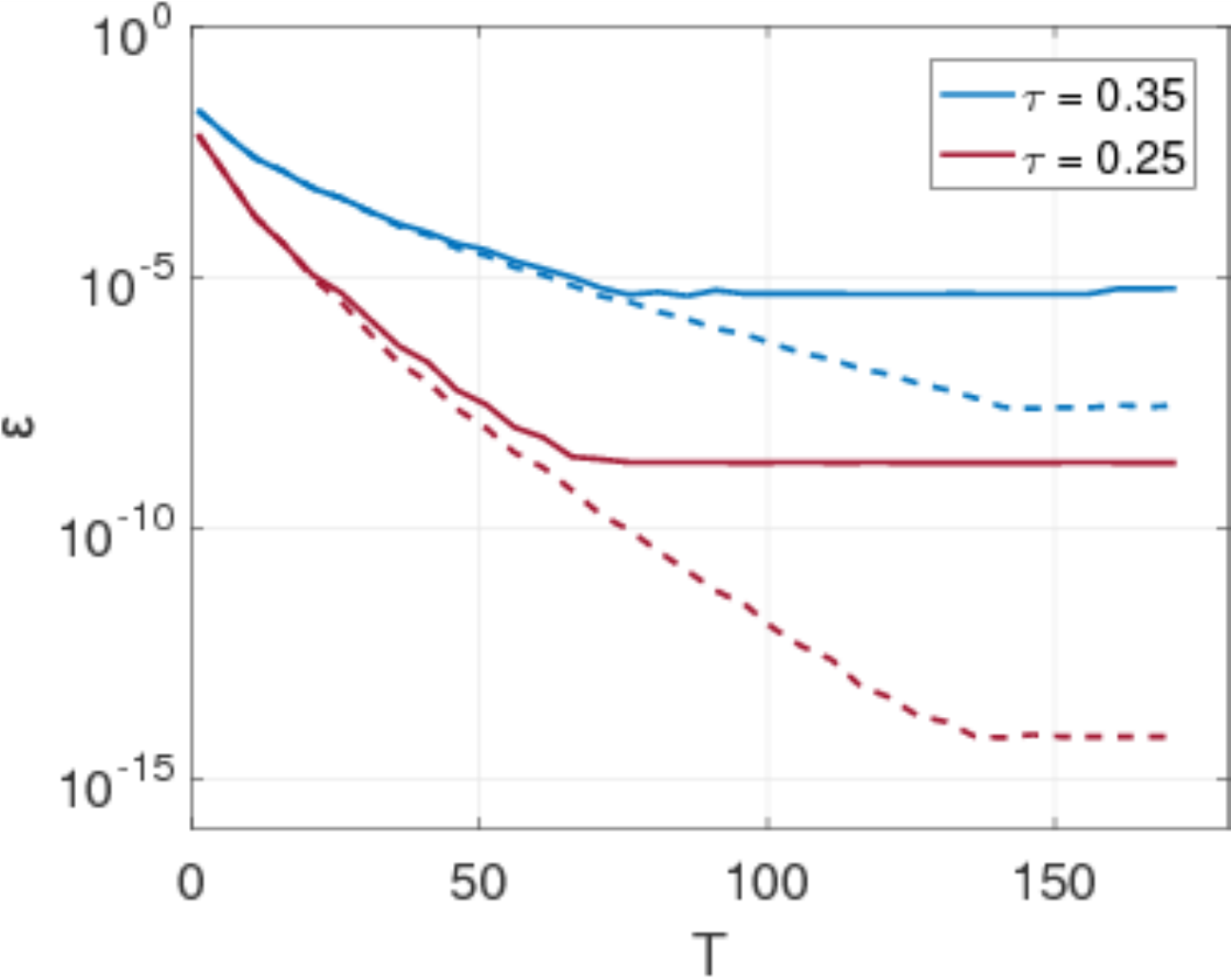}
  		\caption{Approximation error}
  		\label{fig:InvElConvg}
  	\end{subfigure}
  	\begin{subfigure}[d]{0.45\textwidth}
  		\includegraphics[width=\textwidth]{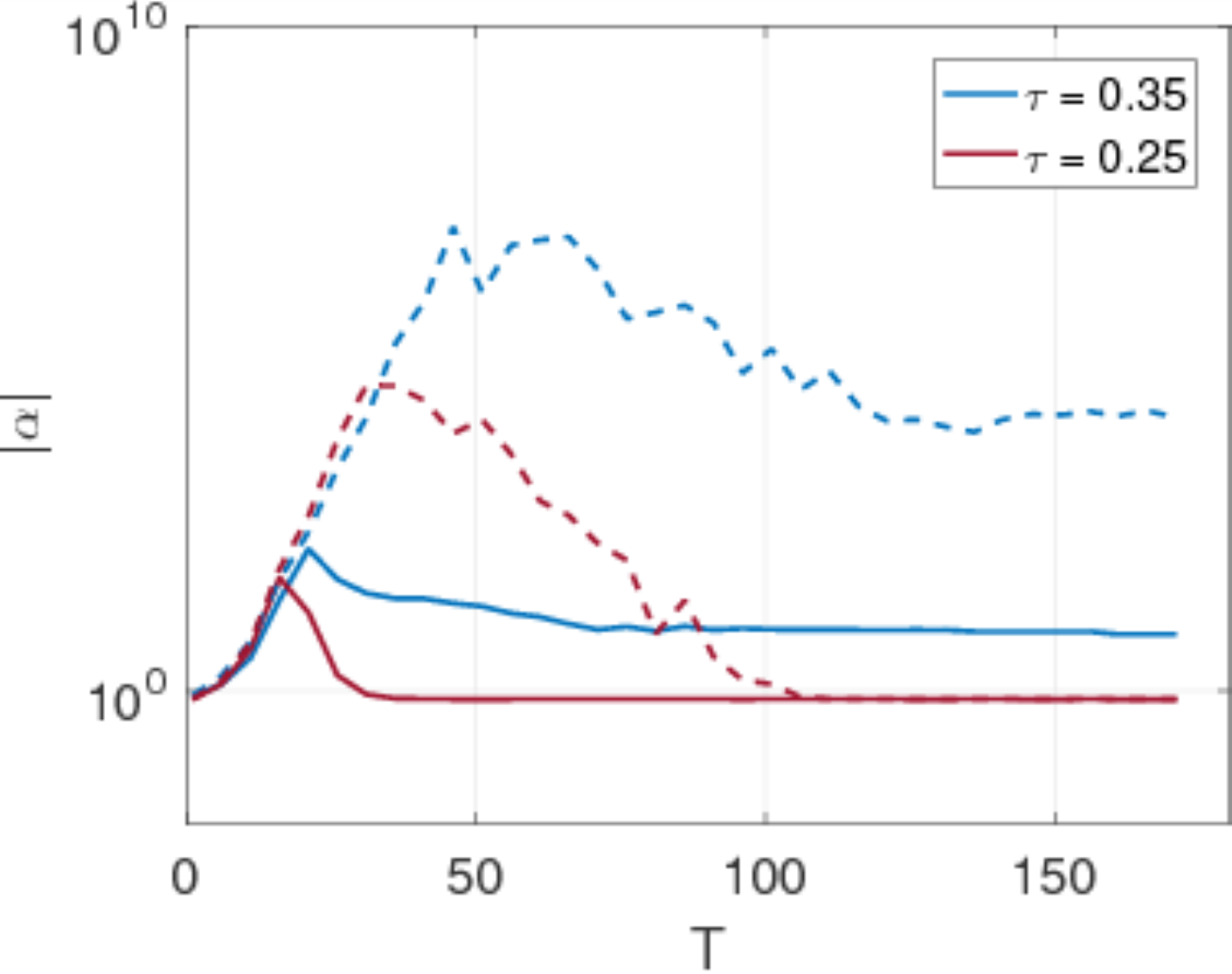}
  		\caption{$2-$norm of the coefficient vector}
  		\label{fig:InvElCoeff}
  	\end{subfigure}      	
  	\caption{WBM results using weighted residual (line) and collocation (dashed line) formulations with constant boundary data for the inverted ellipse.}
  	\label{fig:InvElResults}
  \end{figure}

The corresponding simulation results are shown in Fig.~\ref{fig:InvElResults}. At this stage, the pattern looks familiar. The collocation method achieves higher accuracy. Although the convergence rate is rather slow, and although the corresponding matrices for increasing $T$ become larger and increasingly ill-conditioned, for the first ellipse the method ultimately achieves machine precision accuracy. This does not happen for the second ellipse, where the Schwarz singularity is inside the box, and accuracy is limited to about $10^{-7}$. The limiting size of the coefficient norm is correspondingly larger. The results for weighted residual formulation are similar, but with maximal accuracy of $\sqrt{\epsilon}$ and $10^{-5}$ respectively and, as before, with a smaller limiting norm of the solution vector compared to the collocation approach.

In this final experiment we have increased the oversampling factor to $\gamma=4$. This may be due to the fact that we chose equispaced points in the parameter domain $t \in [0,2\pi]$, yet the corresponding points on the boundary $\partial \Omega$ are unevenly distributed. Still, recall that the unknown constants $\kappa_{M,N}^\epsilon$ and $\lambda_{M,N}^\epsilon$ in \eqref{eq:error_collocation} are known to be bounded only in the limit of $M \to \infty$. As of yet, there is no theoretical guarantee that linear oversampling is sufficient, unless $\gamma$ can be arbitrarily increased when convergence is not seen.

 
\section{Discussion}
\label{sect:discussion}

We focus on the interpretation of the results in this paper, rather than on the mathematical analysis. Though the set of WBM basis functions does not constitute a frame for the relevant function spaces involved, which we do not show in this paper, the generic error bounds \eqref{eq:error_galerkin} and \eqref{eq:error_collocation} do apply.

\subsection{Convergence to accuracy on the order of $\epsilon$ versus $\sqrt{\epsilon}$}

Arguably, the matter of numerical stability is more important than the distinction between very high and even higher precision. Indeed, our choices of the truncation parameter are much larger than is customary in the literature on WBM, since extremely high accuracy is rarely requested in applications. Still, the bound \eqref{eq:error_collocation} for the discrete least squares approximation shows maximal convergence to a value that is proportional to machine precision, and the results in Fig.~\ref{fig:WaveBox} show that this maximal accuracy is obtained using computations in standard double floating point precision only. This is remarkable, in view of the extreme ill-conditioning of the matrices involved as illustrated in Fig.~\ref{fig:WaveBox}(b). Indeed, at some point the matrices are numerically indistinguishable from being singular.

The results for the weighted residual formulation show convergence to roughly $\sqrt{\epsilon}$, as predicted by \eqref{eq:error_galerkin}. We reiterate that, in order to obtain these results, we resorted to a direct solver that minimises the norm $\Vert x \Vert$ of the solution vector, while also minimising the approximation error which corresponds to the residual. Achieving both objectives simultaneously is possible only due to the redundancy in the discretisation: both systems are numerically \emph{underdetermined}, i.e., the number of singular values of $A$ or $\tilde{A}$ larger than a threshold is actually smaller than the number of degrees of freedom, even when oversampling. This redundancy arises in WBM from restricting the bounding box to a smaller domain and it is crucial: it enables the freedom in a solver to minimise $\Vert x \Vert$, while maintaining a small residual. We refer to \cite{huybrechs2016tw674} for a detailed analysis of this effect.

We add that it is also possible to oversample the weighted residual formulation. Indeed, the number of functions of the form \eqref{eq:wavef} is infinite and we have truncated the set using \eqref{eq:truncation}. The formulation that was used in this paper leads to square matrices, because we used an equal number of functions in the expansion \eqref{eq:wbmSol} as in the weighting function in \eqref{eq:weightRes}. However, we could use more functions for the weight, leading to a rectangular system as in the oversampled collocation approach. A similar bound to \eqref{eq:error_collocation} holds in this case \cite{adcock2018frames2}, with maximal accuracy on the order of $\epsilon$. This could be a viable approach for Trefftz methods based on a variational formulation, for which a collocation approach is not easily available.

\subsection{Singularities of the solution}

As mentioned in the introduction, a crucial question is the following: does a solution to the linear system exist, with small residual and with small norm? In the absence of singularities on the boundary of the domain itself, the answer to that question lies in the location of singularities in the extension of the Helmholtz solution from the domain $\Omega$ to the bounding box $S$ (recall Fig.~\ref{fig:WBMdomains}).

The set of functions in the Wave Based Method is the truncation of a complete set in the bounding box, complete in the sense that the solution to any boundary value problem with Neumann data on the four edges of the rectangle can be represented \cite[\S2.5.2]{desmet98phd}. The extension of the solution of the Helmholtz equation on a domain $\Omega$ to a larger domain is unique, at least when the domain is sufficiently smooth such that it admits a parameterisation that is analytic. Note that the expansion \eqref{eq:wbmSol} satisfies the Helmholtz equation not only on $\Omega$, but on the whole of $S$. Therefore, expansion \eqref{eq:wbmSol} necessarily approximates not only the true solution on $\Omega$, but also its extension to $S$. A sufficient condition for the existence of an approximate expansion with bounded norm coefficients is that this extension is bounded and free of singularities on $S$.

In practice, the situation is slightly more involved. Though the extension of the Helmholtz problem is mathematically unique, there is some `wiggle room', for lack of a better terminology. A weaker condition than the above, but also sufficient, is that an extension exists whose Dirichlet or Neumann trace on $\Omega$ agrees \emph{to within an $\epsilon$ tolerance} with the given boundary condition. Finally, it is sufficient for an approximation to such an extension to exist for finite $N$, corresponding to the chosen value of the truncation parameter $T$, and the extension may vary with $N$ (or equivalently with $T$).

It seems that the \emph{wiggle room} increases quantitatively if the singularity is further away from the domain $\Omega$. An important qualitative difference is whether the singularity of the Helmholtz extension lies within the box $S$ or outside.

These considerations were illustrated in Fig.~\ref{fig:CircleHankel}. We solve an interior problem, for which we know that the unique extension has a point source singularity outside $\Omega$ and we vary the location of that singularity.  Convergence to high accuracy is shown in Fig.~\ref{fig:CircleHankel}(a) when the point source is outside the box. The minimal achievable error deteriorates as the point source moves inside the box and closer to the domain $\Omega$. Correspondingly, the size $\Vert x \Vert$ of the solution vector grows as shown in Fig.~\ref{fig:CircleHankel}(b). The maximal accuracy of the approximation problem is negatively correlated with the distance to the nearest singularity.

The fact that problems exist for which the solver does not converge to $\epsilon$ accuracy indicates numerically that the WBM set of functions is not a frame. Indeed, if it were, then convergence to maximal precision should be seen for sufficiently large $N$ for any boundary condition. Mathematically, the existence of such approximations is precisely guaranteed by the frame property \eqref{eq:frame}.

If nearby singularities arise in a practical computation, the solution is to enrich the approximation space, i.e., to add other solutions of the Helmholtz equation to the expansion. The goal should be to build an approximation space in which the exact solution can be approximated with small norm coefficients.

\subsection{Non-convex domains}

The Wave Based Method does not apply to all non-convex domains, but it applies to some. Exterior extensions of solutions to the Helmholtz equation almost always develop singularities at certain points, determined by the geometry of the domain. If these points are outside the bounding box, then convergence is not at stake. It they are inside the box then, as for convex domains with exterior singularities, there is a maximal achievable accuracy that deteriorates with decreasing distance of the singularity. The underlying cause is the same as in the convex case: convergence is limited by the growth of the norm of the solution vector, i.e., of the size of the coefficients in expansion \eqref{eq:wbmSol}.

Whether exterior singularities present a problem in practice can perhaps be determined a-priori based on physical insights. Unfortunately, the computation of their location for non-convex obstacles is not straightforward, and for this reason it remains advisable to restrict WBM to convex obstacles. Still, if a singularity can be expected inside the box, a possible solution may be to enrich the approximation space. For the method of fundamental solutions, the solution proposed in \cite{barnett2008stability} is to position the charge points in between the domain $\Omega$ and the singularities. It seems plausible to assume that adding such singular solutions to the approximation space of WBM would improve convergence.

Alternatively, the issues associated with non-convex domains and solutions with nearby singularities may be treated using h-refinement (subdividing the computational domain), rather than p-refinement (increasing the number of degrees of freedom on a single domain). Indeed this is an approach taken by several Trefftz methods \cite{hiptmair2016trefftz}. The interplay between h- or hp-refinement on the one hand, and frames and redundancy on the other hand, remains to be investigated.

\section{Conclusions}
\label{sect:conclusions}

Ill-conditioning of the linear systems arising in Trefftz methods is often perceived as a cause of concern. We show in this paper that, for at least one Trefftz method, this need not be so: numerical convergence can be guaranteed for a range of problems, in spite of potentially extreme ill-conditioning. We hasten to add that this is certainly not a general observation, and in particular it is limited only to approximation problems involving an approximation space that is, in some sense that can be made precise, redundant. For any other type of problems, ill-conditioning is, and always will be, a major potential source of loss of accuracy in numerical computations.

We show that accurate solutions are found, if they exist in the span of the chosen basis functions, and if their expansion coefficients are not too large. In turn, this means that the limitations of a method can be uncovered by studying when this is not the case. With the experiments of this paper, we have focused on the presence of singularities. There may be other reasons why the exact solution can not be well approximated with small coefficients, and a prominent difficulty arises in applications involving evanescent waves.

The beneficial aspects of oversampled collocation approaches have been advocated in literature before, e.g. \cite{barnett2008stability}. Enriching approximation spaces to capture properties of a solution (resulting in smaller coefficients in the corresponding approximate expansion) have also been explored, for example by adding corner singularities \cite{deckers2012corners}. The main contribution of this paper is the increased level of (mathematical) confidence with which one can make statements about convergence, accuracy and numerical stability of a Trefftz method in the presence of numerical ill-conditioned linear systems.

Through experiments, we could show applicability of WBM beyond convex domains, a limitation that is commonly accepted in literature, as long as the non-convex domain is not too non-convex. We showed improved accuracy of oversampled collocation. Still, the scope of this paper was limited and the contents led to several topics of ongoing and future research. In particular, we did not examine the influence of corner singularities in the solution. We did not analyse the oversampling factor $\gamma$ (recall that $M = \gamma N$ with $\gamma > 1$). The generic error bounds we used do not indicate the convergence rate of the solution. We have focused solely on p-refinement. Finally, we have studied only a very restricted subset of the many available settings in which the Wave Based Method has been successfully applied.

\section*{Acknowledgements}

The authors greatly appreciate discussions on the topic of this paper with Ben Adcock, Alex Barnett, Timo Betcke, Vincent Copp{\'e}, Elke Deckers, Wim Desmet, Andrew Gibbs, David Hewett, Ralf Hiptmair, Stijn Jonckheere, Roel Matthysen, Andrea Moiola and Marcus Webb. The first author was supported by FWO-Flanders projects G.0641.11 and G.A004.14, as well as by KU Leuven project C14/15/055.

\bibliography{references}

\end{document}